\documentclass[11pt,a4paper]{article}
\usepackage[utf8]{inputenc}
\usepackage{amsmath,amssymb,amsfonts}
\usepackage[margin=3cm]{geometry}
\usepackage{xcolor}
\usepackage{algorithm}
\usepackage{algorithmicx}
\usepackage{algpseudocode}
\usepackage{enumitem}
\usepackage{booktabs}
\usepackage{cite}
\usepackage[hyperfootnotes=false]{hyperref}

\allowdisplaybreaks[2]

\hypersetup{
colorlinks=true,
linkcolor=blue,
citecolor=darkgreen,
urlcolor=blue,
pdftitle={Proximal random reshuffling with inexact gradients for nonconvex composite optimization},
pdfauthor={Boris S. Mordukhovich and Pedro P\'erez-Aros}}
\definecolor{darkgreen}{RGB}{0,150,0}

\newcommand{\Csta}{\mathtt{a}}
\newcommand{\Cstb}{\mathtt{b}}
\newcommand{\R}{\mathbb R}
\def\oR{\overline{\R}}
\newcommand{\N}{\mathbb N}
\newcommand{\PLKname}{Polyak--{\L}ojasiewicz--Kurdyka}
\newcommand{\prox}[2]{\mathtt{prox}_{#1 #2}}
\DeclareMathOperator{\dom}{dom}
\DeclareMathOperator{\dist}{dist}

\usepackage{amsthm}
\newtheorem{theorem}{Theorem}[section]
\newtheorem{lemma}[theorem]{Lemma}
\newtheorem{definition}[theorem]{Definition}
\newtheorem{assumption}[theorem]{Assumption}
\theoremstyle{remark}
\newtheorem{remark}[theorem]{Remark}

\title{Proximal random reshuffling with inexact gradients for nonconvex composite optimization}

\author{Boris S.\ Mordukhovich\thanks{Department of Mathematics and Institute for AI and Data Science,
                             Wayne State University,
                             Detroit, MI 48202, \textsc{USA}.
                                 Email:~\href{mailto:aa1086@wayne.edu}
                                 {aa1086@wayne.edu}}
        \and
        Pedro P\'erez-Aros\thanks{Departamento de Ingenier\'ia Matem\'atica and Centro de Modelamiento Matem\'atico (CNRS IRL2807), Universidad de Chile, Santiago, \textsc{Chile}.
        Email:~\href{mailto:pperez@dim.uchile.cl}{pperez@dim.uchile.cl}}
 }

\makeatletter
\renewcommand{\maketitle}{
  \begin{center}
    {\LARGE\bfseries\@title\par}
    \vskip 2em
    {\large
      \lineskip .5em
      \begin{tabular}[t]{c}
        \@author
      \end{tabular}\par}
  \end{center}
  \@thanks
}
\makeatother

\begin{document}

\maketitle

\begin{abstract}
We design and justify a proximal inexact random reshuffling method for minimizing
sums of the finite average of functions with Lipschitz continuous gradients and lower semicontinuous regularizers that may be
nonconvex. Each epoch visits every component once, using inexact gradients with a prescribed tolerance, and ends with a single proximal step. Our analysis is pathwise: it holds for every sequence of permutations and does not depend on how the reshuffling is performed. We show that when the stepsizes and
error tolerances satisfy suitable summability conditions, the objective values converge and every accumulation point is stationary. Furthermore, under the Polyak--{\L}ojasiewicz--Kurdyka condition at a single accumulation point and sufficient decrease of stepsizes, the
entire sequence has finite length and converges to that point. For polynomial stepsizes and power-type desingularizing functions, we derive convergence rates for objective values and iterates. Finally,  faster objective-value rates are derived when the limit is a local minimizer.
\end{abstract}


\section{Introduction}\label{intro}
{\em Finite-sum objective functions} arise in large-scale learning,
statistical estimation, stochastic programming, and other applications in which the objective aggregates many individual loss functions; see. e.g.,
\cite{BottouCurtisNocedal2018} for an overview and further references. In many of such applications, the sum is complemented by a {\em regularization term}, which encodes constraints and structural features of the sought solution being often {\em nonsmooth} and {\em nonconvex}. Problems combining a large sum with a possibly nonconvex regularizer form therefore a mathematically natural
setting of wide applicability; see the books
\cite{Beck2017FOM,sra2011optimization,Lan2020FirstOrderStochastic} and the bibliographies therein. This motivates us to study the following composite problem:
\begin{align}\label{Problem:P}\tag{P}
\min_{x\in\mathbb R^n}\ \varphi(x):=f(x)+h(x),
\qquad f(x):=\frac{1}{\mathtt N}\sum_{i=1}^{\mathtt N}f_i(x),
\end{align}
where each $f_i:\mathbb R^n\to\mathbb R$ has  Lipschitz continuous
gradient and $h:\mathbb R^n\to\mathbb R\cup\{\infty\}$ is a
lower semicontinuous (l.s.c.) regularizer. Convexity of $h$ is not assumed while $h$
is required to be prox-bounded and continuous relative to its domain, and $\varphi$ is bounded from below; see Assumption~\ref{Assumpt01}. Stationarity for \eqref{Problem:P} is
understood in the sense of the basic (limiting/Mordukhovich) subdifferential, i.e., as the $M$-{\em stationarity} $0\in\partial\varphi(\bar x)$.

When $\mathtt N$ is large, evaluating $\nabla f$ at every iteration is expensive, and incremental methods update iterates by using one component at a time. In the convex setting, incremental subgradient methods and incremental proximal methods, in which proximal steps are applied to individual components, are analyzed in \cite{NedicBertsekas2002} and \cite{Bertsekas2011IncrementalProximal},
respectively. {\em Stochastic gradient} methods (SGMs) in the line of \cite{RobbinsMonro1951} sample the components with replacement, whereas
{\em random reshuffling} (RR) visits every component exactly once per epoch according to  permutations typically redrawn at each epoch. For strongly convex smooth finite sums, the analysis in \cite{MR4214476} shows that RR with diminishing stepsizes
converges faster than $1/t$, while the $1/t$ rate cannot be improved in general for with-replacement sampling. In the smooth nonconvex case, \cite{SecciaCoppolaLiuzziPalagi2025} proves stationarity of accumulation points of incremental and reshuffling gradient schemes under gradient Lipschitz continuity, with convergence being enforced by a
watchdog rule and a sporadically activated linesearch rather than by a growth property of the objective.

For nonconvex composite problems, the standard tool for proving convergence of the entire sequence of iterates, rather than merely
subsequential stationarity, is an inequality bounding the objective gap at a point by a function of the (sub)gradient norm at that point. Such an inequality was introduced in \cite{POLYAK1963864}: for a smooth
function $f$, the condition $\|\nabla f(x)\|^2\geq2\mu\,(f(x)-\inf f)$
with $\mu>0$ was shown there to yield linear convergence of the gradient method without any convexity assumptions. Independently,
\cite{lojasiewicz65} established, for real-analytic functions, the gradient inequality $|f(x)-f(\bar x)|^{\theta}\le C\|\nabla f(x)\|$, $\theta\in[\frac12,1[$, near a critical point $\bar x$. Later on,
\cite{kurdyka} replaced the power function by a general desingularizing function $\psi$ for functions definable in o-minimal structures obtaining the estimate $\|\nabla(\psi\circ f)(x)\|\ge 1$ and, as a consequence, the finite length of gradient trajectories. An nonsmooth extension via the Clarke subdifferential was introduced \cite{bolte2007}. Conditions of this type have been used extensively to prove convergence of iterates generated by descent optimization methods; see, e.g., \cite{attouch13,Bento2025,fund} with the bibliographies therein. In these works, the main results show that sequences of iterates satisfying sufficient-decrease and relative-error conditions on the subgradients converge to stationary points as soon as the objectives satisfy the aforementioned gradient-type error bounds at one of theiraccumulation points. Following \cite{Bento2025,fund}, we label the latter properties as the {\em Polyak--{\L}ojasiewicz--Kurdyka} (PLK) conditions (see Definition~\ref{def_PKL}), a name that records the three contributions above and, in particular, the pioneering role of~\cite{POLYAK1963864} in the analysis of descent methods.

{\em Gradient errors}, a second feature of the method studied here, arise from approximate derivative computations, inexact inner solvers, or perturbed search directions. For smooth convex problems, \cite{DevolderGlineurNesterov2014} introduces an inexact-oracle model
and quantifies the accumulation of oracle errors in first-order methods. In the smooth nonconvex setting, \cite{zbMATH07960830}
analyzes inexact reduced gradient methods with controlled gradient errors, establishes stationarity of accumulation points and---under PLK conditions----convergence of the entire sequence with rates. Related developments in \cite{khanh2024fundamental} provide a fundamental convergence analysis of sharpness-aware minimization and its variants
by interpreting these algorithms within an inexact gradient descent framework and including methods with nonmonotone objective values. The finite-length convergence analysis in \cite{khanh2024fundamental} relies on an auxiliary decreasing quantity and a quasi-additivity-type property of the desingularizing function. Although this property holds for power-type desingularizing functions, it constitutes an additional structural assumption in the general
setting; see \cite[Assumption~2.3]{khanh2024fundamental}. For composite problems with a convex regularizer, \cite{BollapragadaGupta2026} studies proximal gradient methods with adaptively controlled gradient accuracy and derives convergence with complexity guarantees. These works do  not address the within-epoch accumulation of gradient errors arising in reshuffling methods.

Reshuffling methods are generally {\em nonmonotone}: the objective value may not decrease from one epoch to the next. Consequently, the conventional descent frameworks developed in \cite{attouch13,Bento2025} and many other publications cannot be applied. The analysis in \cite{MR4605214} overcomes this difficulty for smooth
RR with exact gradients and diminishing stepsizes. Assuming that each component function is bounded from below with Lipschitz gradient, \cite{MR4605214}  establishes
convergence of the entire sequence to a stationary point under an appropriate
PLK condition, together with convergence rates depending on the PLK exponent. Its  convergence analysis relies on a quasi-additivity-type property of the desingularizing function. In \cite{JoszLaiLi2026}, a proximal RR is studied for locally Lipschitzian functions with convex regularizers without assuming either coercivity or the existence of accumulation points. The conclusions in \cite{JoszLaiLi2026}, obtained through a tracking lemma relating the iterates to trajectories of a continuous-time dynamical system, concern near-approximate stationarity rather than convergence of the entire sequence. Proximal (sub)gradient methods in which neither summand is required to be either convex or smooth are analyzed in \cite{Solodov2025ProximalSubgradient} under Clarke-regularity-type assumptions, with establishing stationarity of accumulation points as the main conclusion.

The combination of hypotheses and conclusions considered in our paper is therefore not covered by the previous results. {\em{First}}, the regularizer is convex
in \cite{JoszLaiLi2026,BollapragadaGupta2026} and is absent in \cite{MR4605214}. In contrast, for an l.s.c. and possibly nonconvex regularizer $h$ as considered here,  the stationarity residual must be measured through proximal optimality conditions and the subdifferential, and a one-sided PLK condition has to be reconciled with iterates whose
objective values may cross the level $\varphi(\bar x)$. {\em Second}, inexact
component gradients and the epoch error inherent in reshuffling must be controlled simultaneously. {\em Third}, the componentwise lower boundedness used in \cite{MR4605214} to control the within-epoch drift excludes situations in which only the composite objective is bounded below. We replace it by the aggregate gradient-growth condition in Assumption~\ref{Assumpt01}\ref{Assumpt01:item3}, which is implied by componentwise lower boundedness; see Remark~\ref{rem:growth-generalizes}. {\em Fourth}, the quasi-additivity property of the desingularizing
function, imposed in \cite{MR4605214,khanh2024fundamental},
is not guaranteed by the PLK conditions and may be difficult to verify for general nonsmooth composite optimization problems. Indeed, the existence results for desingularizing functions
in \cite{kurdyka,bolte2007} are not constructive and do not necessarily provide an explicit representation of such
functions in terms of the data. This contrasts with the classical gradient
inequalities in \cite{lojasiewicz65,POLYAK1963864},
which admit power-type desingularizing functions for which quasi-additivity can be readily established. Consequently, although quasi-additivity holds for important classes of desingularizing functions, it cannot generally be inferred from the PLK condition alone.\vspace*{0.03in}

The {\em main contributions} of our paper are the following.

\begin{enumerate}
\item We introduce the novel {\em proximal inexact random reshuffling} (PIRR)
method, Algorithm~\ref{alg:1}, which combines random reshuffling of the components, inexact component gradients with tolerance
$\varepsilon_t$, and a single proximal step per epoch with respect to an l.s.c.\  and possibly nonconvex regularizer $h$.

\item An {\em approximate version} of the {\em descent inequality} for PIRR, with explicit coefficients, is established in Theorem~\ref{lemma:descent}. Convergence of the {\em objective values} and summability of the {\em weighted increments} and {\em stationarity residual} are proved in Theorem~\ref{prop:01}, which are instrumental for deriving {\em stationarity of every accumulation point} in Theorem~\ref{thm:cluster-stationarity}.

\item Convergence of the {\em entire sequence} of iterates under the basic {\em PLK condition} at a single accumulation point is established in Theorem~\ref{Theo_main_convergencePLK}. No quasi-additivity-type property of the desingularizing function and/or any additional assumptions are imposed. The {\em key
ingredient} that makes it possible is the rescaling result of Lemma~\ref{lem:noqap-rescaling}. Remarks~\ref{rem:value-attainment}
and~\ref{rem:smooth-no-step-ratio} discuss the continuity hypothesis on $h$ and the smooth case where  $h\equiv0$.

\item {\em Convergence rates} for the objective values and the iterates under polynomial stepsizes and the exponent PLK condition are given in Theorem~\ref{main_rates}, and {\em faster objective-value rates} when the limit is a local minimizer are presented in Theorem~\ref{thm:improved-local-rates}.  Table~\ref{tab:optimal-rates} summarizes
the corresponding  optimal choices of the stepsize exponent.
\end{enumerate}\vspace*{0.05in}

The rest of the paper is organized as follows. Section~\ref{Preliminaries} fixes the notation,
collects auxiliary results on real sequences, and introduces PIRR. Section~\ref{Sect:SubsequenceConver} contains the within-epoch discrepancy estimate, the
approximate descent inequality together with the subsequential convergence and stationarity results. Section~\ref{sec:plk} introduces the basic PLK condition, proves the rescaling lemma and the finite-length theorem, and derives the convergence rates with the weakly convex and local minimizer refinements. Section~\ref{sec:7} discusses limitations and future research. Appendices~\ref{AppendixA} and~\ref{AppendixB} contain the
proofs of two auxiliary lemmas on recursive inequalities.

\section{Preliminaries and Proximal Inexact Random Reshuffling Algorithm}\label{Preliminaries}
In this section, we introduce the basic notation and definitions used in the rest of this work. We also recall several preliminary results that will be instrumental in conducting our analysis. Finally, we present our novel numerical method formalized as the \emph{Proximal Inexact Random Reshuffling} algorithm.

\subsection{Basic notation}

In this work, we use conventional notation from optimization and variational analysis (see, e.g., \cite{MR1491362,MR4659163,MR3823783}). Denote by $\R^n$ the $n$-dimensional Euclidean space equipped with the Euclidean norm $\|\cdot\|$, and by $\mathbb{N}:=\{1,2,\ldots\}$ the set of positive integers. Unless otherwise stated, all extended-real-valued functions $f\colon\R^n\to\oR:=]-\infty,\infty]$ are supposed to be {\em proper}, i.e., $\dom f:=\{x\in\R^n\;|\;f(x)<\infty\}\ne\emptyset$.

A differentiable function $f:\R^n\to\R$ is said to be {\em $\mathtt{L}$--smooth} if its gradient is is Lipschitz continuous with modulus $\mathtt{L}\ge 0$ ($\mathtt{L}$--Lipschitz continuous), i.e.,
\[
\|\nabla f(x)-\nabla f(y)\|\le \mathtt{L}\|x-y\|\;\mbox{ for all }\;x,y\in\R^n.
\]
Given a permutation $\sigma$ of $\{1,\ldots,\mathtt{N}\}$ (i.e., a bijection mapping of $\{1,\ldots,\mathtt{N}\}$ onto itself), we use the shorthand
\[
\sigma_i:=\sigma(i), \qquad i=1,\ldots,\mathtt{N}.
\]
For a sequence of nonnegative numbers $(s_t)_{t\ge 1}$ and a function $g:]0,\infty[\to]0,\infty[$, we write $s_t=\mathcal{O}(g(t))$ as $t\to\infty$ if there exist $C>0$ and $t_0\ge 1$ such that
\[
s_t\le C\,g(t)\quad \text{for all }\,t\ge t_0.
\]
Furthermore, we  write $s_t=\mathcal{O}^{-}(t^{-\nu})$ if $s_t=\mathcal{O}(t^{-\nu'})$ for all $\nu'<\nu$.

Let $h:\R^n\to\oR$ be an l.s.c.\ function,  and let $\lambda>0$. The \emph{proximal operator} of $h$ with parameter $\lambda$ is the (possibly set-valued) mapping
\[
\prox{\lambda}{h}(x):=\arg\min_{y\in\R^n}\left\{\, h(y)+\frac{1}{2\lambda}\|y-x\|^2 \right\},
\qquad x\in\R^n.
\]
A function $h:\R^n\to\oR$ is \emph{prox-bounded} if there exists $\lambda>0$ ensuring that\\ $\inf_{y\in\R^n}\left\{\, h(y)+\frac{1}{2\lambda}\|y-x\|^2 \right\}>-\infty$ for every $x\in\R^n$.
The supremum of all such parameters $\lambda$, denoted by $\lambda_h$, is called the \emph{threshold of prox-boundedness} for $h$.

The \emph{basic $($limiting/Mordukhovich$)$ subdifferential} of $\varphi\colon\R^n\to\oR$ at $\bar x\in\dom\varphi$ is
\begin{equation}\label{sub}
\partial \varphi(\bar x)
:=\big\{v\in\R^n\;\big|\;\exists\, x^k\to \bar x,\ \varphi(x^k)\to \varphi(\bar x),\ v^k\to v,\
v^k\in \widehat{\partial}\varphi(x^k),\; k\in\N\big\}
\end{equation}
with $\partial\varphi(\bar x):=\emptyset$ if $\bar x\notin\dom\varphi$, where  $ \widehat{\partial}\varphi( x)$ stands for the \emph{regular $($Fr\'echet$)$ subdifferential}, which is defined by $\widehat{\partial}\varphi(\bar x)
:=\big\{v\in\R^n\;\big|\;
\liminf_{x\to \bar x}
\frac{\varphi(x)-\varphi(\bar x)-\langle v, x-\bar x\rangle}{\|x-\bar x\|}\ge 0\big\}$
for $\bar x \in \dom \varphi$ and  $\widehat{\partial}\varphi(\bar x):=\emptyset$ whenever $\bar x\notin\dom\varphi$.

Given a point $x\in\mathbb{R}^n$ and a set $S\subset\mathbb{R}^n$, we denote the {\em distance} from $x$ to $S$ by  $\dist(x;S):=\inf\left\{\|x-y\|\;|\;y\in S\right\}$. We say that $\bar x$ is a {\em stationary point} for problem \eqref{Problem:P} provided that $0 \in \partial \varphi(\bar x)$ via the basic subdifferential \eqref{sub}.

\subsection{Preliminary results}
To analyze convergence of our algorithm, we begin with recalling the standard \emph{descent lemma} for smooth functions; see, e.g., \cite[Theorem~3.6]{MR4605214}. This inequality will be used repeatedly in our analysis.

\begin{lemma}\label{descent}
Let $h\colon\mathbb{R}^n \to \mathbb{R}$ be a differentiable function that is $\mathtt{L}$-smooth. Then for any $x, y \in \mathbb{R}^n$, we have the inequality
\begin{equation*}
h(y) \leq h(x) + \langle \nabla h(x), y - x \rangle + \frac{\mathtt{L}}{2} \|x - y\|^2.
\end{equation*}
\end{lemma}
We also recall a deterministic {\em Robbins--Siegmund sequence lemma} taken from \cite{MR343355}.

\begin{lemma}\label{Prop:Robbins–Siegmund}
Let $\{v_t\}$, $\{\beta_t\}$, $\{\xi_t\}$, and $\{\eta_t\}$ be sequences of nonnegative  numbers. Suppose that for all $t \geq 1$ it holds that
\begin{align*}
v_{t+1}\leq (1 +\beta_t) v_t + \xi_t - \eta_t.
\end{align*}
Then the sequence $\{v_t\}$ converges and we have the series convergence $\sum_{t=1}^{\infty} \eta_t < \infty$  provided that $\sum_{t=1}^{\infty} \beta_t < \infty$ and $\sum_{t=1}^{\infty}\xi_t <\infty$.
\end{lemma}

The following lemma established in \cite[Lemma~3.1]{zbMATH07960830} plays a crucial role in the subsequent convergence analysis.

\begin{lemma}\label{LemmaKham2023}
Let $\{x^t\}_{t\in\N}, \{d^t\}_{t\in\N} \subset \R^n$ satisfy the condition $$
\sum_{t=1}^{\infty}\|x^{t+1}-x^t\|\cdot\|d^t\|<\infty.
$$
Assume that $\bar x$ is an accumulation point of $\{x^t\}_{t\in\N}$ and that zero is an accumulation point of $\{d^t\}_{t\in\N}$. Then there exists a subsequence $\{t_j\}_{j\in\N}\subset\N$ such that
\[
x^{t_j}\to \bar x
\quad\text{and}\quad
d^{t_j}\to 0
\quad\text{as }\;j\to\infty.
\]
\end{lemma}

Next we present three auxiliary results on real sequences. The first two lemmas below provide convergence rate estimates for sequences satisfying recursive inequalities. These results are technical variants of
\cite[Chapter~2.2, Lemmas~4 and 5]{zbMATH07325960}. We give their elaborated proofs in Appendices~\ref{AppendixA} and \ref{AppendixB}, respectively, for completeness and reader's convenience. Lemma~\ref{lem:rate_sequence} is the tool used to extract the convergence rates of Theorem~\ref{main_rates} from the scalar
recursions obtained in its proof, whereas
Lemma~\ref{lem:local-rate-scalar} is tailored to the recursion that arises when the accumulation point is a local minimizer and is used in the
proof of Theorem~\ref{thm:improved-local-rates}.

\begin{lemma}\label{lem:rate_sequence}
Let $\{s_t\}$ be a nonnegative sequence satisfying the inequality
\begin{equation}\label{eq:implicit-scalar-recursion}
s_t+\frac{\lambda}{(t+b)^\gamma}s_t^q
\leq s_{t-1}+\frac{D}{(t+b)^{\gamma+q\mu}}
\end{equation}
for all large $t$, where $\lambda>0$, $D,b\geq0$, $\mu>0$, $q\geq1$, and
$\gamma\in]0,1]$. Then we have
\[
s_t=
\begin{cases}
\mathcal O(t^{-\mu}),
&0<\gamma<1,\quad q=1,\\
\mathcal O\left(
t^{-\min\{\mu,\frac{1-\gamma}{q-1}\}}
\right),
& 0<\gamma<1,\quad q>1,\\
\mathcal O(t^{-\mu}),
& \gamma=q=1,\quad \lambda>\mu,\\
\mathcal O(t^{-\lambda}),
& \gamma=q=1,\quad \lambda<\mu,\\
\mathcal O\left(
[\ln(t+1)]^{-\frac1{q-1}}
\right),
& \gamma=1,\quad q>1.
\end{cases}
\]
\end{lemma}

\begin{lemma}\label{lem:local-rate-scalar}
Let $\{s_t\}$ be a nonnegative sequence satisfying
\begin{equation}\label{eq:local-scalar-recursion}
s_t+\frac{\lambda}{(t+b)^\gamma}s_t^q
\le s_{t-1}+\frac{D}{(t+b)^{3\gamma}}
\end{equation}
for all large $t$, where $\lambda>0$, $D,b\geq0$,
$0<q\leq1$, $1/2<\gamma\leq1$, and $(q,\gamma)\neq(1,1)$. Then we have the rate $s_t=\mathcal O(t^{-2\gamma/q})$.
\end{lemma}

The final lemma is an estimation of power series, which is derived by using the classical integral comparison. This lemma provides tail estimates for polynomial stepsizes that are needed in the proof of Theorem~\ref{Theo_main_convergencePLK} to
bound accumulated perturbations and in the proof of Theorem~\ref{main_rates} to
verify the imposed assumptions.

\begin{lemma}\label{Lemma_3.7rate}
Given   $a>0$, $b\geq0$, and $\gamma\in]1/3,1]$, suppose that
\begin{equation}\label{eq:polynomial-schedule}
\alpha_t=\frac{a}{(t+b)^\gamma}\quad\mbox{and}\quad
\varepsilon_t=\mathcal O((t+b)^{-\gamma}).
\end{equation}
Then we have $	\sum_{j=t}^{\infty}\bigl(\alpha_j^3+\alpha_j\varepsilon_j^2\bigr)
=\mathcal O((t+b)^{1-3\gamma})$. Moreover, it holds for $\theta\in]0,1[$ that
\begin{equation}\label{eq:weighted-power-tail}
\sum_{j=t}^{\infty}\alpha_j
\left[\sum_{i=j}^{\infty}(\alpha_i^3+\alpha_i\varepsilon_i^2)\right]^\theta=\mathcal O\bigl((t+b)^{-(\gamma+\theta(3\gamma-1)-1)}\bigr)
\end{equation}
provided that $\gamma+\theta(3\gamma-1)>1$.
\end{lemma}\vspace*{-0.05in}

\subsection{The Proximal Inexact Random Reshuffling Algorithm}
Now are ready to design the proximal inexact random reshuffling algorithm (Algorithm~\ref{alg:1}) for solving problem~\eqref{Problem:P} by the following procedure.
\begin{algorithm}[htbp]
\begin{algorithmic}[1]
\Require{$x^0 \in \dom h$}
\For{$t = 1, 2, \ldots$}
\State Choose stepsize $\alpha_t > 0$
\State Generate a (random) permutation $\sigma^t$ of $\{1, \ldots, \mathtt{N}\}$
\State Set $\tilde{x}_0^t := x^{t-1}$
\For{$i = 1, \ldots, \mathtt{N}$}
\State Compute $g_i^t$ such that $\| g_i^t - \nabla f_{ \sigma_i^t}(\tilde{x}_{i-1}^t) \| \leq \varepsilon_t$
\State Update $\tilde{x}_i^t = \tilde{x}_{i-1}^t - \alpha_t g_i^t$
\EndFor
\State Choose $x^t \in \prox{\alpha_t \mathtt{N} }{h}\left( \tilde{x}_{\mathtt{N}}^t \right)$
\EndFor
\end{algorithmic}
\caption{Proximal Inexact Random Reshuffling}
\label{alg:1}
\end{algorithm}

Note that in each epoch $t$ of Algorithm~\ref{alg:1}, the inner iterates
$\tilde x_i^t$ start from $x^{t-1}$ and visit the components in the order given by $\sigma^t$ while using {\em inexact gradients} $g_i^t$ with error at most $\varepsilon_t$. The epoch is closed by a single proximal step with parameter $\alpha_t\mathtt N$. The permutation $\sigma^t$ may be either fixed
for all epochs, or drawn at random as in random reshuffling, where it is sampled from a prescribed distribution. Observe that our analysis does not depend on how $\sigma^t$ is generated. Furthermore, two consequences of this structure drive the analysis in the next
section. Since the inner iterates move away from $x^{t-1}$, the averaged direction $\frac1{\mathtt N}\sum_{i=1}^{\mathtt N}g_i^t$ differs from $\nabla f(x^{t-1})$ even when $\varepsilon_t=0$. This epoch error and
the tolerance $\varepsilon_t$ are estimated together in Lemma~\ref{estiDist:gradient}. Moreover, the optimality condition in the
proximal step relates the stationarity residual
$\dist(0;\partial\varphi(x^t))$ to the increment $x^t-x^{t-1}$, which is
the content of Lemma~\ref{estiDist}.

\section{Subsequence Convergence Analysis for PIRR Algorithm}\label{Sect:SubsequenceConver}

In this section, we study convergence properties of the iterative sequence generated by Algorithm~\ref{alg:1}. We begin with stating the main assumptions on the problem data in \eqref{Problem:P}.

\begin{assumption} \label{Assumpt01}
Let $f=\mathtt N^{-1}\sum_{i=1}^{\mathtt N}f_i$ and
$\varphi=f+h$. Assume that:
\begin{enumerate}[label=\textup{(\roman*)}]
\item
Each $f_i:\R^n\to\R$ has a globally
$\mathtt L_i$-Lipschitz continuous gradient.

\item The regularizer $h:\R^n\to\oR$ is l.s.c.  and prox-bounded. We also suppose that
$h$ is continuous relative to $\dom h$ and that the objective is bounded from below, i.e.,
$\varphi(x)\geq\underline\varphi$ for some $\underline\varphi\in\R$ and all $x\in\dom h$.

\item \label{Assumpt01:item3} There are constants $\Csta,\Cstb\geq0$ such that
\begin{equation}\label{eq:aggregate-growth}
\sum_{i=1}^{\mathtt N}\|\nabla f_i(x)\|^2
\leq \Csta+\Cstb\bigl(\varphi(x)-\underline\varphi\bigr)\;\mbox{ for all }\;x\in\dom h.
\end{equation}

\item The proximal subproblems defining the
algorithm have minimizers, which holds eventually whenever the positive
parameters $\mathtt N\alpha_t$ lie below the prox-boundedness threshold.
\end{enumerate}
\end{assumption}

\begin{remark}\label{rem:growth-generalizes}
It is worth mentioning that the constant $\underline\varphi$ in  Assumption~\ref{Assumpt01} may not coincide with $\inf\varphi$. Moreover, this assumption is satisfied whenever all the component
functions $f_i$ and the regularizer $h$ are bounded from below (under the remaining smoothness and regularity hypotheses). Indeed, applying Lemma~\ref{descent} to each $f_i$ at
$y=x-\mathtt L_i^{-1}\nabla f_i(x)$ with $\mathtt L_i>0$ yields the estimate
\[
\|\nabla f_i(x)\|^2
\leq 2\mathtt L_i\bigl(f_i(x)-\inf f_i\bigr).
\]
Summing up the latter inequalities brings us for all $x\in\dom h$ to
\[
\sum_{i=1}^{\mathtt N}\|\nabla f_i(x)\|^2
\leq 2\mathtt N\max_i\mathtt L_i
\bigl(\varphi(x)-\underline\varphi\bigr)\;\text{ with } 	\underline\varphi
:=\frac1{\mathtt N}\sum_{i=1}^{\mathtt N}\inf f_i+\inf h.
\]
Note that the converse inequality may not hold. For example, setting
$f_1(x):=x$, $f_2(x):=-x$, and $h(x):=x^2$ on $\R$ leads us to
$\varphi(x)=x^2$ and $\sum_i|f_i'(x)|^2=2$, although neither component is bounded from below. Finally, if $\dom h$ is bounded, the continuity of the gradients $\nabla f_i$ implies that their squared norms are uniformly bounded on $\dom h$. Thus \eqref{eq:aggregate-growth} holds with $\Cstb=0$ and some $\Csta\geq0$.
\end{remark}

\subsection{Approximate descent property and stationarity of accumulation points} In this subsection, we establish important results for Algorithm~\ref{alg:1} concerning an approximate descent inequality, value convergence, and stationarity of accumulation points of iterates. For all the results of this subsection, suppose that {\em Assumption~{\rm\ref{Assumpt01}} holds}. We first derive an estimate of the {\em stationarity residual}.

\begin{lemma}\label{estiDist}
Let the sequence of iterates $\{x^t\}_{t\geq0}$ be generated by Algorithm~{\rm\ref{alg:1}}. Then for any $t\geq1$ and $\nu>0$, we have the estimate
\[
\begin{aligned}
\dist^2\!\left(0;\partial\varphi(x^t)\right)
&\le(1+\nu)\left(\frac1{\mathtt N\alpha_t}
+\frac1{\mathtt N}\sum_{i=1}^{\mathtt N}\mathtt L_i\right)^2
\|x^t-x^{t-1}\|^2\\
&\quad+(1+\nu^{-1})
\left\|\nabla f(x^{t-1})
-\frac1{\mathtt N}\sum_{i=1}^{\mathtt N}g_i^t
\right\|^2.
\end{aligned}
\]
\end{lemma}
\begin{proof}
By the necessary optimality condition for the subproblem in the proximal step defining $x^t$ and the subdifferential sum rule from \cite{MR3823783}, we deduce the inclusion
\[
0 \in \partial h(x^t)+\frac{1}{\alpha_t \mathtt N}(x^t-x^{t-1})+\frac{1}{\mathtt N}\sum_{i=1}^{\mathtt N} g_i^t.
\]
Adding and subtracting $\nabla f(x^{t-1})$ and using $\partial \varphi(x^t)=\nabla f(x^t)+\partial h(x^t)$ yield
\[
-\left(\frac{1}{\alpha_t \mathtt N}(x^t-x^{t-1})+\frac{1}{\mathtt N}\sum_{i=1}^{\mathtt N} g_i^t-\nabla f(x^{t-1})\right)
+\bigl(\nabla f(x^t)-\nabla f(x^{t-1})\bigr)
\in \partial \varphi(x^t),
\]
which therefore brings us to the inequality
\begin{align*}
\dist\!\left(0;\partial \varphi(x^t)\right)
&\le\frac{1}{\alpha_t\mathtt N}\|x^t-x^{t-1}\|+\left\|\nabla f(x^{t-1})-
\frac1{\mathtt N}\sum_{i=1}^{\mathtt N}g_i^t\right\|\\
&\quad+\|\nabla f(x^t)-\nabla f(x^{t-1})\|.
\end{align*}
Remembering that each $\nabla f_i$ is $\mathtt L_i$-Lipschitz tells us that
\begin{align*}
\|\nabla f(x^t)-\nabla f(x^{t-1})\|\le  \frac{1}{\mathtt N}\left(\sum_{i=1}^{ \mathtt N} \mathtt L_i\right)\|x^t-x^{t-1}\|.
\end{align*}
Then we complete the proof by using
$(a+b)^2\le (1+\nu)a^2+(1+\nu^{-1})b^2$.
\end{proof}

The next lemma estimates the discrepancy between the {\em averaged inexact gradients} and the full one in Algorithm~\ref{alg:1}.

\begin{lemma}\label{estiDist:gradient}
Let $\{x^t\}_{t\geq0}$ be generated by Algorithm~{\rm\ref{alg:1}}. If $\alpha_t\to0$, then there exists $t_0\in\N$ such that for every $t\geq t_0$, we have
\begin{equation}\label{eq:direct-epoch-discrepancy}
\begin{aligned}
\left\|
\frac1{\mathtt N}\sum_{i=1}^{\mathtt N}g_i^t
-\nabla f(x^{t-1})\right\|^2 \le 4\alpha_t^2\left(\sum_{i=1}^{\mathtt N}\mathtt L_i^2\right)
\left[\Csta+\Cstb
\bigl(\varphi(x^{t-1})-\underline\varphi\bigr)\right]
+4\varepsilon_t^2.
\end{aligned}
\end{equation}
\end{lemma}
\begin{proof}
For $i=1,\ldots,\mathtt N$, define $w_i^t:=g_i^t-\nabla f_{\sigma_i^t}(x^{t-1})$.
The inner updates satisfy $\tilde x_{i-1}^t-x^{t-1}=-\alpha_t\sum_{j=1}^{i-1}g_j^t$.
The gradient error bounds and Lipschitz continuity of the component gradients combined with the Cauchy--Schwarz inequality give us
\[
\begin{aligned}
\|w_i^t\|&\le\mathtt L_{\sigma_i^t}
\|\tilde x_{i-1}^t-x^{t-1}\|+\varepsilon_t\leq
\alpha_t\mathtt L_{\sigma_i^t}
\left\|\sum_{j=1}^{i-1}g_j^t\right\|+\varepsilon_t\leq
\alpha_t \mathtt L_{\sigma_i^t}
\sum_{j=1}^{\mathtt N}\|g_j^t\|
+\varepsilon_t\\
&\le\alpha_t \mathtt L_{\sigma_i^t}
\sum_{j=1}^{\mathtt N}\|w _j^t\|  + \alpha_t \mathtt L_{\sigma_i^t}
\sum_{j=1}^{\mathtt N}\|\nabla f_j (x^{t-1})\|
+\varepsilon_t.
\end{aligned}
\]
In this way, we arrive at the estimates
\[
\begin{aligned}
\|w_i^t\|^2 &\le 3 \left(\alpha^2_t \mathtt{ L}^2_{\sigma_i^t} \left[
\sum_{j=1}^{\mathtt N}\|w _j^t\| \right]^2 + \alpha^2_t \mathtt{ L}^2_{\sigma_i^t}
\left[\sum_{j=1}^{\mathtt N}\|\nabla f_j (x^{t-1})\|\right]^2
+\varepsilon^2_t \right) \\
&\le 3 \left(\alpha^2_t \mathtt{ L}^2_{\sigma_i^t}  \mathtt N
\sum_{j=1}^{\mathtt N}\|w _j^t\|^2  + \alpha^2_t \mathtt{ L}^2_{\sigma_i^t}
\mathtt N 	\sum_{j=1}^{\mathtt N}\|\nabla f_j (x^{t-1})\|^2+\varepsilon^2_t \right)\\
&\le 3\left(\alpha^2_t \mathtt{ L}^2_{\sigma_i^t}  \mathtt N
\sum_{j=1}^{\mathtt N}\|w _j^t\|^2  + \alpha^2_t\mathtt{ L}^2_{\sigma_i^t}
\mathtt N\left(\Csta+\Cstb\bigl(\varphi(x^{t-1})-\underline\varphi\bigr) \right)
+\varepsilon^2_t \right),
\end{aligned}
\]
where the last inequality uses \eqref{eq:aggregate-growth}. Recalling that $\sigma^t$ is a permutation, we get
\[
\begin{aligned}
\sum_{i=1}^{\mathtt N}\|w_i^t\|^2
&\leq 3 \alpha^2_t  \mathtt{N}   \sum_{i=1}^{\mathtt N}\mathtt L_i^2\left( \sum_{j=1}^{\mathtt N}\|w_j^t\|^2 \right) +3 \alpha^2_t  \mathtt{N}   \sum_{i=1}^{\mathtt N}\mathtt L_i^2 \left(\Csta+\Cstb\bigl(\varphi(x^{t-1})-\underline\varphi\bigr) \right) + 3\mathtt{N} \varepsilon_t^2.
\end{aligned}
\]
It follows that $3 \alpha^2_t  \mathtt{N}   \sum_{i=1}^{\mathtt N}\mathtt L_i^2
\leq\frac14$ when $t\in\N$ is sufficiently large. Therefore,
\[
\begin{aligned}
\sum_{i=1}^{\mathtt N}\|w_i^t\|^2
\leq 4 \alpha^2_t\mathtt{N}   \sum_{i=1}^{\mathtt N}\mathtt L_i^2 \left(   \sum_{j=1}^{\mathtt N}\|\nabla f_j (x^{t-1})\|^2 \right)+ 4\mathtt{N}\varepsilon_t^2.
\end{aligned}
\]
Remembering that $\frac1{\mathtt N}\sum_{i=1}^{\mathtt N}g_i^t-\nabla f(x^{t-1})
=\frac1{\mathtt N}\sum_{i=1}^{\mathtt N}w_i^t$ yields
\begin{equation}
\begin{aligned}
\left\|
\frac1{\mathtt N}\sum_{i=1}^{\mathtt N}g_i^t
-\nabla f(x^{t-1})\right\|^2&\leq
4\alpha_t^2\left(\sum_{i=1}^{\mathtt N}\mathtt L_i^2\right)\left[\Csta+\Cstb
\bigl(\varphi(x^{t-1})-\underline\varphi\bigr)\right]
+4\varepsilon_t^2,
\end{aligned}
\end{equation}
which thus completes the proof of the lemma.
\end{proof}

Now we obtain the main results of this subsection. First we establish a crucial {\em approximate descent property} for the sequence of iterates generated by Algorithm~\ref{alg:1}.

\begin{theorem}
\label{lemma:descent} Let $\{x^t\}_{t\geq0}$ be generated by Algorithm~{\rm\ref{alg:1}}, and let $\alpha_t\to 0$ as $t\to\infty$. Then for any fixed $\delta\in]0,1/2[$, there exist $t_\delta\in\mathbb N$ and $C_\delta>0$ such that
\begin{equation}\label{EQ01}
\begin{aligned}
\varphi(x^t)\hspace{-0.05cm}-\hspace{-0.05cm}\underline\varphi
&\leq(1\hspace{-0.05cm}+\hspace{-0.05cm}\mathfrak a_t)
\bigl(\varphi(x^{t-1})\hspace{-0.05cm}-\hspace{-0.05cm}\underline\varphi\bigr)
\hspace{-0.05cm}-\hspace{-0.05cm}\mathfrak b_t\dist^2\bigl(0;\partial\varphi(x^t)\bigr)\hspace{-0.05cm}-\hspace{-0.05cm}\mathfrak c_t\|x^t-x^{t-1}\|^2+\mathfrak{d}_t
\end{aligned}
\end{equation}
whenever $t\geq t_\delta$ with the explicitly calculated constants
\begin{equation}\label{defConstants}
\begin{aligned}
\mathfrak a_t
&=C_\delta \alpha_t^3 , \quad
\mathfrak b_t=\left(\frac{1}{2}-\delta \right)\mathtt N\alpha_t,
\quad\mathfrak c_t=\frac{\delta}{2\mathtt N\alpha_t},\quad\mathfrak{d}_t
=C_\delta\left[\alpha_t^3
+\alpha_t\varepsilon_t^2\right].
\end{aligned}
\end{equation}
\end{theorem}
\begin{proof}
Denote $\mathtt{L}: = \frac1{\mathtt N}\sum_{i=1}^{\mathtt N}\mathtt L_i$  and deduce from the definition of the proximal operator the following inequality
\[
\begin{aligned}
h(x^t)&+\frac{\|x^t-x^{t-1}\|^2}{2\mathtt N\alpha_t}+\left\langle
\frac1{\mathtt N}\sum_{i=1}^{\mathtt N}g_i^t,
x^t-x^{t-1}\right\rangle
\leq h(x^{t-1}).
\end{aligned}
\]
Combining this inequality with Lemma~\ref{descent} for $f$ ensures that
\begin{equation}\label{eq:descent-prox-basic}
\begin{aligned}
&\varphi(x^{t-1})-\varphi(x^t)
\ge\left(\frac1{2\mathtt N\alpha_t}-\frac{\mathtt L}{2}\right)\|x^t-x^{t-1}\|^2\\
&+\left\langle\frac1{\mathtt N}\sum_{i=1}^{\mathtt N}g_i^t
-\nabla f(x^{t-1}),x^t-x^{t-1}
\right\rangle\\
&\ge\frac{1}{\mathtt N \alpha_t}	\left(\frac{1-\delta/2}{2}-\frac{\mathtt{N} \alpha_t \mathtt L}{2}\right)\|x^t-x^{t-1}\|^2\quad- \frac{\mathtt N\alpha_t}{\delta}\left\|
\frac1{\mathtt N}\sum_{i=1}^{\mathtt N}g_i^t
-\nabla f(x^{t-1})\right\|^2.
\end{aligned}
\end{equation}
Then we get from Lemma~\ref{estiDist} with $\eta>0$ that
\begin{equation}\label{eq:descent-prox-basic_2}
\begin{aligned}
\mathfrak{b}_t\dist^2(0;\partial\varphi(x^t))
&\le\left(\frac{1}{2} -\delta\right)(1+ \eta)\left(1+\alpha_t \mathtt{N}\mathtt L
\right)^2\frac{\|x^t-x^{t-1}\|^2}{\mathtt N\alpha_t}\\
&\quad+\mathtt N \alpha_t	\left(1+\eta^{-1}\right)\left\|
\frac1{\mathtt N}\sum_{i=1}^{\mathtt N}g_i^t
-\nabla f(x^{t-1})\right\|^2.
\end{aligned}
\end{equation}
Employing \eqref{eq:descent-prox-basic} and \eqref{eq:descent-prox-basic_2} and noticing that
\[\lim\limits_{ t\to \infty ,  \eta \to 0  }\left(\frac{1-\delta/2}{2}-\frac{ \mathtt{N}\alpha_t \mathtt L}{2}\right)-\left(\frac{1}{2} -\delta\right)(1+ \eta)\left(
1+\mathtt{N}\alpha_t\mathtt{L}
\right)^2 =  \frac{3}{4}\delta
\]
ensure the existence of constants $C_\delta>0$ and $t_\delta$ such that
\begin{equation}\label{eq:sharp-two-term}
\begin{aligned}
\varphi(x^{t-1})-\varphi(x^t)&\geq
\left(\frac12-\delta\right)\mathtt N\alpha_t
\dist^2(0;\partial\varphi(x^t)) +
\frac{\delta}{2\mathtt N\alpha_t}
\|x^t-x^{t-1}\|^2\\
&\quad-C_\delta\mathtt N\alpha_t
\left\|\frac1{\mathtt N}\sum_{i=1}^{\mathtt N}g_i^t-\nabla f(x^{t-1})\right\|^2\;\mbox{ for all }\;t\geq t_\delta.
\end{aligned}
\end{equation}
Finally, we substitute \eqref{eq:direct-epoch-discrepancy} into \eqref{eq:sharp-two-term} and rearrange the terms to verify \eqref{EQ01} by increasing $C_\delta$ if needed to match
with the coefficients in \eqref{defConstants}.
\end{proof}

\begin{remark}\label{rem:weakly-convex-descent} Here we discuss two issued related to Theorem~\ref{lemma:descent}.

{\bf(i)} It has been well realized that {\em decent properties} (in particular, {\em sufficient decrease}) {\em fail} to hold for SGD and RR methods. However, such properties play a crucial role in establishing convergence with rates in the variety of deterministic algorithms; see \cite{attouch13,Bento2025,fund}. The novel {\em approximate} version of the decent property proposed in Theorem~\ref{lemma:descent} allows us to proceed with the convergence analysis of PIRR in what follows; cf.\ also some counterparts for exact algorithms in \cite{khanh2024fundamental,MR4605214}.

{\bf(ii)} Suppose that in the setting of
Theorem~\ref{lemma:descent}, $h$ is {\em globally $\rho$-weakly convex} with some $\rho\ge 0$, i.e., the function
$x\mapsto h(x)+\frac{\rho}{2}\|x\|^2$
is convex. Then \eqref{EQ01} holds with the {\em improved coefficient} $\mathfrak b_t=(1-\delta)\mathtt N\alpha_t,$ while $\mathfrak c_t$ and $\mathfrak a_t,\mathfrak{d}_t$ retain the forms in \eqref{defConstants} after enlarging $t_\delta$ and $C_\delta$ if needed. Indeed, the proximal optimality condition gives us
\[
-\frac{x^t-x^{t-1}}{\mathtt N\alpha_t}
-\frac1{\mathtt N}\sum_{i=1}^{\mathtt N}g_i^t
\in\partial h(x^t).
\]
Applying the subgradient inequality for weakly convex functions at the comparison point $x^{t-1}$ and combining the resulting
estimate with Lemma~\ref{descent} for $f$, yield
\[
\begin{aligned}
\varphi(x^{t-1})-\varphi(x^t)
&\geq
\left(\frac1{\mathtt N\alpha_t}
-\frac{\mathtt L+\rho}{2}\right)
\|x^t-x^{t-1}\|^2\hspace{-0.1cm}+\hspace{-0.1cm}
\left\langle
\frac1{\mathtt N}\sum_{i=1}^{\mathtt N}g_i^t
\hspace{-0.1cm}	-\hspace{-0.1cm}\nabla f(x^{t-1}),
x^t\hspace{-0.1cm}-\hspace{-0.1cm}x^{t-1}
\right\rangle.
\end{aligned}
\]
Therefore, the remaining proof follows the same arguments as in Theorem~\ref{lemma:descent}.
\end{remark}

The next theorem imposes {\em summability conditions} on the stepsizes and gradient errors to derive a descent-type inequality for the objective values along the iterates and to establish the corresponding {\em value convergence} properties.

\begin{theorem}\label{prop:01}
Let $\sum_t\alpha_t^3<\infty$, $\sum_t\alpha_t\varepsilon_t^2<\infty$ and fix
$\delta\in]0,1/2[$. Then the value sequence $\{\varphi(x^t)\}$ in Algorithm~{\rm\ref{alg:1}} converges to a finite number, and we have
\begin{equation}\label{eq00:prop:01}
\sum_{t=1}^\infty\frac{\|x^t-x^{t-1}\|^2}{\alpha_t}<\infty,
\qquad
\sum_{t=1}^\infty\alpha_t
\dist^2(0;\partial\varphi(x^t))<\infty.
\end{equation}
Moreover, there exists a positive number $n_0>0$ ensuring the estimate
\begin{equation}\label{ineq:nondecreasing}
\begin{aligned}
\varphi(x^t)&\leq\varphi(x^{t-1})
-\mathfrak b_t\dist^2(0;\partial\varphi(x^t))
-\mathfrak c_t\|x^t-x^{t-1}\|^2+n_0\bigl(\alpha_t^3+\alpha_t\varepsilon_t^2\bigr)
\end{aligned}
\end{equation}
for all large $t\in\N$, where $\mathfrak b_t$ and $\mathfrak c_t$ are taken from \eqref{defConstants}.
\end{theorem}
\begin{proof}
Define $v_t:=\varphi(x^t)-\underline\varphi\ge 0$ and deduce from \eqref{EQ01} that
\[
v_t\leq(1+\mathfrak a_t)v_{t-1}+\mathfrak{d}_t
-\mathfrak b_t\dist^2(0;\partial\varphi(x^t))
-\mathfrak c_t\|x^t-x^{t-1}\|^2
\]
for all large $t$. The coefficients $\mathfrak a_t$ and $\mathfrak{d}_t$ are
nonnegative and summable due to $\alpha_t\to0$,
$\sum_t\alpha_t^3<\infty$, and
$\sum_t\alpha_t\varepsilon_t^2<\infty$. Therefore, Lemma~\ref{Prop:Robbins–Siegmund} yields the convergence of $v_t$ and summability of the two
nonnegative descent terms in \eqref{eq00:prop:01}.  Furthermore, since $v_t$ is bounded, the perturbation
$\mathfrak a_tv_{t-1}+\mathfrak{d}_t$ is bounded from above by a fixed multiple of $\alpha_t^3+\alpha_t\varepsilon_t^2$.
Substitution this into \eqref{EQ01} gives us
\eqref{ineq:nondecreasing}.
\end{proof}

Now we are in a position to establish {\em subsequential convergence} of the iterative sequence $\{x^t\}$ in Algorithm~\ref{alg:1}
to a {\em stationary point} of problem \eqref{Problem:P}.

\begin{theorem}\label{thm:cluster-stationarity} In the setting of Algorithm~{\rm\ref{alg:1}}, suppose that
\[
\sum_{t=1}^{\infty}\alpha_t=\infty,
\quad\sum_{t=1}^{\infty}\alpha_t^3<\infty,
\quad\sum_{t=1}^{\infty}\alpha_t\varepsilon_t^2<
\infty.
\]
Then every accumulation point $\bar x$ of $\{x^t\}$ satisfies the stationary condition $0\in\partial\varphi(\bar x)$.
\end{theorem}
\begin{proof}
Choose a subgradient $d^t\in\partial\varphi(x^t)$ with $\|d^t\|=\dist(0;\partial\varphi(x^t))$, which is possible by the proximal optimality condition and closedness of the basic
subdifferential at $x^t$. It follows from Theorem~\ref{prop:01} and Young's inequality that
\[
\sum_{t=1}^{\infty}
\|x^{t+1}-x^t\|\cdot\|d^{t+1}\|
\leq\frac12\sum_{t=1}^{\infty}
\left(\frac{\|x^{t+1}-x^t\|^2}{\alpha_{t+1}}
+\alpha_{t+1}\|d^{t+1}\|^2\right)<\infty.
\]
Since $\sum_t\alpha_t\|d^t\|^2<\infty$ while
$\sum_t\alpha_t=\infty$, zero is an accumulation point of $\{d^t\}$. Fix an arbitrary accumulation point $\bar x$ of $\{x^t\}$. Applying Lemma~\ref{LemmaKham2023} to $\{x^t\}$ and the shifted sequence $\{d^{t+1}\}$ allows us to find subsequences of $\{x^t\}$ and $\{d^{t+1}\}$ such that $x^{t_j}\to\bar x$ and $d^{t_j+1}\to 0 $. We have by \eqref{eq00:prop:01} that $\|x^{t+1}-x^t\|\to 0$, and so
$x^{t_j+1}\to\bar x$. The imposed relative continuity property of $h$ yields $\varphi(x^{t_j+1})\to\varphi(\bar x)$. Passing to the limit as $j\to\infty$ and using the closedness of the basic subdifferential graph brings us to $0\in\partial\varphi(\bar x)$ and thus completes the proof of the theorem.
\end{proof}

\section{Convergence under PLK}\label{sec:plk}
In this section we investigate convergence of the entire sequence $\{x^t\}$ generated by Algorithm~\ref{alg:1} under an appropriate version of the \PLKname\ condition introduced as follows; cf.\ \cite{attouch13,Bento2025,fund}.

\begin{definition}\label{def_PKL} We say that an l.s.c.\ function $\varphi:\R^n\to\oR$ satisfies the $($basic$)$ PLK condition at $\bar x\in\dom\varphi$ if there exist a
neighborhood $U$ of $\bar x$, a positive constant $\eta$, and a continuous concave function $\psi:[0,\eta[\to[0,\infty[$ such that $\psi(0)=0$,
$\psi\in C^1]0,\eta[$ with $\psi'>0$ on $]0,\eta[$, and
\begin{equation}\label{PKL_2}
\psi'(\varphi(x)-\varphi(\bar x))
\dist(0;\partial\varphi(x))\ge 1
\end{equation}
whenever $x\in U$ and $0<\varphi(x)-\varphi(\bar x)<\eta$.
\end{definition}
Semialgebraic objectives and the like provide familiar examples of desingularizing supporting PLK \cite{attouch13,bolte2007}. Note that we do not impose the {\em quasi-additivity} estimate
\begin{equation}\label{quasi}
C(\psi'(s+r))^{-1} \leq (\psi'(s))^{-1} + (\psi'(r))^{-1}
\end{equation}
with some constant $C>0$ as, e.g., in \cite{MR4605214,zbMATH07960830,khanh2024fundamental}. The next lemma shows that a rescaled desingularizing function satisfies an appropriate substitute of \eqref{quasi}, which suffices to control the error terms in \eqref{ineq:nondecreasing}.

\begin{lemma} \label{lem:noqap-rescaling}
Let $\psi:[0,\eta[\to[0,\infty[$ be continuous and concave with $\psi(0)=0$, $\psi\in C^1(]0,\eta[)$, and $\psi'>0$ on $]0,\eta[$. Fix $k>0$, put $c:=k+1$, and pick $\tau\in]0,\eta[$, $u\in[0,k\tau]$, and $v\in\mathbb R$ such that
$v+u\in]0,\eta[$. Then we have the inequality
\begin{equation}\label{eq:psqadi}
\left[\psi'\left(\frac{u+v}{c}\right)\right]^{-1}
\le
\begin{cases}
\left[\psi'(v)\right]^{-1}+
\left[\psi'(\tau)\right]^{-1}
& \text{if }\;v\in]0,\eta[,\\
\left[\psi'(\tau)\right]^{-1}
& \text{if }\;v\le 0.
\end{cases}
\end{equation}
\end{lemma}
\begin{proof}
Define $a:=v+u$ and observe that it follows from the concavity and differentiability of $\psi$ on $]0,\eta[$ that $s\mapsto \left[\psi'(s)\right]^{-1}$ is nondecreasing on $]0,\eta[$. From the fact that $a\in]0,\eta[$ and $c\geq1$, we get $a/c\in]0,\eta[$, and so $\psi'(a/c)$ is well defined.

First we examine the case where $v\in]0,\eta[$. If $ a\leq c\tau$, then $\left[\psi'\left(\frac{a}{c}\right)\right]^{-1}\le\left[\psi'(\tau)\right]^{-1}$, which verifies \eqref{eq:psqadi}. If $ a > c\tau$, then using $u\leq k\tau$ and $c=k+1$ tells us that
$u\le k\tau<k a c^{-1}= k (k+1)^{-1} a =a-ac^{-1}$, and so $v=a-u>\frac{a}{c}$. By the monotonicity of
$s\mapsto[\psi'(s)]^{-1}$, we get  $\left[\psi'\left(\frac{a}{c}\right)\right]^{-1}\le
\left[\psi'(v)\right]^{-1}$, which also verifies \eqref{eq:psqadi}. If $v\le 0$, it follows that $a=v+u\leq u\leq k\tau$. Consequently, $ a c^{-1} \leq k (k+1)^{-1} \tau \leq \tau$, which yields $\left[\psi'\left(\frac{a}{c}\right)\right]^{-1}\le
\left[\psi'(\tau)\right]^{-1}$ and thus completes the proof.
\end{proof}

Now we are ready to establish the principal result of the paper about the finite-length convergence of Algorithm~\ref{alg:1} under the PLK condition.

\begin{theorem}\label{Theo_main_convergencePLK}
Let $\varphi=f+h$ satisfy Assumption~{\rm\ref{Assumpt01}}, and let
$\{x^t\}$ be generated by Algorithm~{\rm\ref{alg:1}} with the starting point $x^0\in\operatorname{dom}h$.
Assume that the positive stepsizes and nonnegative error tolerances satisfy
\begin{equation}\label{eq01:prop:stat2}
\sum_{t=1}^{\infty}\alpha_t=\infty,\qquad
\sum_{t=1}^{\infty}\alpha_t^3<\infty,\qquad
\sum_{t=1}^{\infty}\alpha_t\varepsilon_t^2<\infty,
\end{equation}
and that, for some $\rho>0$, it holds that
\begin{equation}\label{eq:noqap-step-ratio}
\alpha_{t+1}\leq\rho\alpha_t
\quad\text{ whenever $t$ is sufficiently large}.
\end{equation}
Let $\bar x$ be an accumulation point of $\{x^t\}$ at which $\varphi$ fulfills the PLK condition with a desingularizing function $\psi$. Suppose also that
\begin{equation}\label{desing_condi}
\sum_{t=p_0}^{\infty}
{\alpha_t}\left[ {\psi'\!\left(\sum_{j=t}^{\infty}\bigl(\alpha_j^3+\alpha_j\varepsilon_j^2\bigr)\right)}\right]^{-1}<\infty
\end{equation}
for some $p_0\in\N$ sufficiently large.	Then we have the relationships
\begin{equation}\label{eq:noqap-conclusions}
\sum_{t=1}^{\infty}\|x^t-x^{t-1}\|<\infty,\qquad
x^t\longrightarrow\bar x,\qquad \varphi(x^t) \to \varphi(\bar x), \quad
0\in\partial\varphi(\bar x).
\end{equation}
\end{theorem}
\begin{proof} Fix $\delta\in]0,1/2[$, and let $U$ and $\eta>0$ be the neighborhood and parameter in the PLK condition for $\varphi$ from Definition~\ref{def_PKL} at $\bar x$. Theorem~\ref{prop:01} gives us
\begin{equation}\label{eq:noqap-raw-descent}
\varphi(x^{t+1})
\leq \varphi(x^t)-\mathfrak c_{t+1}\|x^{t+1}-x^t\|^2
-\mathfrak b_{t+1}\dist^2\bigl(0;\partial\varphi(x^{t+1})\bigr)
+\mathtt W_t
\end{equation}
for large $t\in\N$, where $	\mathtt W_t
=n_0\bigl(\alpha_{t+1}^3+\alpha_{t+1}\varepsilon_{t+1}^2\bigr)$, $n_0\geq1$, and such that
\begin{equation}\label{eq:noqap-eventual-coefficients}
\mathfrak c_t\geq\frac{\mathfrak r}{\alpha_t},\qquad
\mathfrak b_t\geq\mathfrak r\alpha_t
\qquad\text{eventually}
\end{equation}
with some constant $\mathfrak r>0$. Define $\mathtt u_t$ below and estimate it by
\begin{equation*}
\mathtt u_t:=\sum_{j=t}^{\infty}\mathtt W_j
=n_0\sum_{j=t+1}^{\infty}
\bigl(\alpha_j^3+\alpha_j\varepsilon_j^2\bigr)
\le n_0\sum_{j=t}^{\infty}
\bigl(\alpha_j^3+\alpha_j\varepsilon_j^2\bigr),
\end{equation*}
where the quantities are finite, and hence $\mathtt u_t\to 0$ as $t\to\infty$. Since $\mathtt u_t-\mathtt u_{t+1}=\mathtt W_t$, we deduce from \eqref{eq:noqap-raw-descent} and \eqref{eq:noqap-eventual-coefficients} that
\begin{equation}\label{eq:noqap-corrected-descent}
\bigl(\varphi(x^t)+\mathtt u_t\bigr)
-\bigl(\varphi(x^{t+1})+\mathtt u_{t+1}\bigr) \geq\mathfrak r\Bigl(\frac{\|x^{t+1}-x^t\|^2}{\alpha_{t+1}}+\alpha_{t+1}\dist^2\bigl(0;\partial\varphi(x^{t+1})\bigr)\Bigr).
\end{equation}
The objective $\varphi$ is bounded from below, and hence the tail of $\bigl(\varphi(x^t)+\mathtt u_t\bigr)$ is nonincreasing and also bounded from below. This tells us that
\begin{equation}\label{eq:noqap-value-limit}
\varphi(x^t)+\mathtt u_t\downarrow\ell,\qquad
\varphi(x^t)\longrightarrow\ell\;\mbox{ as }\;t\to\infty
\end{equation}
for some finite number $\ell$. It readily follows from
\eqref{eq:noqap-corrected-descent} that
\begin{equation*}
\sum_{t=1}^{\infty}\frac{\|x^t-x^{t-1}\|^2}{\alpha_t}<\infty,\qquad
\sum_{t=1}^{\infty}\alpha_t \dist^2\bigl(0;\partial\varphi(x^t)\bigr)<\infty.
\end{equation*}
Since $\alpha_t\to0$ by \eqref{eq01:prop:stat2}, we get the convergence
\begin{equation}\label{eq:noqap-small-Gamma}
\begin{gathered}
\|x^t-x^{t-1}\|\to 0,\qquad
\alpha_t\dist\bigl(0;\partial\varphi(x^t)\bigr)\to 0,\\
\Gamma_t:=\|x^t-x^{t-1}\|+\alpha_t\dist\bigl(0;\partial\varphi(x^t)\bigr)\longrightarrow 0\;\mbox{ as }\;t\to\infty.
\end{gathered}
\end{equation}
Choosing a subsequence $x^{t_j}\to\bar x$ and employing the lower semicontinuity of $\varphi$ give us $\varphi(\bar x)\leq\ell<\infty$, and so $\bar x\in\operatorname{dom}h$. The imposed continuity of $f$ and continuity of $h$ relative to its domain ensure the convergence
$\varphi(x^{t_j})\longrightarrow\varphi(\bar x)$ as $j\to\infty$.
Consequently, we have $\ell=\varphi(\bar x)$, which confirms the value convergence in \eqref{eq:noqap-conclusions}. Therefore, it follows from \eqref{eq:noqap-value-limit} that $\varphi(x^t)-\varphi(\bar x)+\mathtt u_t\ge 0$ and $\varphi(x^t)-\varphi(\bar x)+\mathtt u_t\downarrow 0$.
Moreover, the existence of some $t_0 \in \N$ with
$\varphi(x^{t_0})-\varphi(\bar x)+\mathtt u_{t_0} =0$ implies by \eqref{eq:noqap-corrected-descent} that $x^{t} = x^{t_0}$ for all $t\geq t_0$, and thus $\bar x $ is a stationary point of problem \eqref{Problem:P}. Therefore, for the rest of the proof, we can assume that $\varphi(x^t)-\varphi(\bar x)+\mathtt u_t>0$. Let us set
\begin{align*}
c:=n_0+1, \; \Psi(s):=c\,\psi(s/c),\;
\mathfrak{s}_t:=\varphi(x^t)-\varphi(\bar x)+\mathtt u_t,\;\widehat\Delta_t :=\Psi\!\left(\mathfrak{s}_t\right).
\end{align*}
Since $0<\mathfrak{s}_t \downarrow0$, the definition of $\widehat\Delta_t$ is valid for all large $t$, and
we have $\widehat\Delta_t\downarrow0$.\vspace*{0.03in}

We split the remainder of the proof into the following four claims.\vspace*{0.03in}

{\textbf{Claim~1:}} {\em Let $t\in \mathbb{N}$ such that $\mathfrak{s}_t \in ]0,\eta[$ and $x^t\in U$. Then we have}
\begin{equation}\label{eq:noqap-key-estimate}
\left[{\Psi'\!\left(\mathfrak{s}_t\right)} \right]^{-1}\leq\dist\bigl(0;\partial\varphi(x^t)\bigr)+\left[\psi'\!\left(\sum_{j=t}^{\infty}\bigl(\alpha_j^3+\alpha_j\varepsilon_j^2\bigr)\right)\right]^{-1}.
\end{equation}
To verify this, we use Lemma~\ref{lem:noqap-rescaling} by choosing there
\[
\begin{gathered}
v=\varphi(x^t)-\varphi(\bar x),\qquad
u=\mathtt u_t,\qquad
k=n_0,\quad\tau=\sum_{j=t}^{\infty}\bigl(\alpha_j^3+\alpha_j\varepsilon_j^2\bigr).
\end{gathered}
\]
Observe that $\mathtt u_t$ has the upper and lower bounds
\[
0\leq\mathtt u_t\leq n_0\sum_{j=t}^{\infty}
\bigl(\alpha_j^3+\alpha_j\varepsilon_j^2\bigr).
\]
If $v=\varphi(x^t)-\varphi(\bar x)>0$, then
\[
0<v\leq\varphi(x^t)-\varphi(\bar x)+\mathtt u_t<\eta,
\]
and the PLK inequality \eqref{PKL_2} gives us
\[
\frac{1}{\psi'(v)}
\le\dist\bigl(0;\partial\varphi(x^t)\bigr).
\]
Therefore, Lemma~\ref{lem:noqap-rescaling}  yields the claimed estimate \eqref{eq:noqap-key-estimate}.\vspace*{0.03in}

{\textbf{Claim~2:}} {\em Let $t\in \mathbb{N}$ such that $\mathfrak{s}_t \in ]0,\eta[$ and $x^t\in U$. Then we have}
\begin{equation}\label{eq:noqap-Gamma-recursion}
\Gamma_{t+1}
\leq\left( \Gamma_t-\Gamma_{t+1}\right)
+\frac{\alpha_t}{\psi'\!\left(\sum_{j=t}^{\infty}\bigl(\alpha_j^3+\alpha_j\varepsilon_j^2\bigr)\right)}
+\frac{2\rho}{\mathfrak r}
(\widehat\Delta_t-\widehat\Delta_{t+1}).
\end{equation}
Indeed, the concavity of $\Psi$ tells us that
\[
\widehat\Delta_t-\widehat\Delta_{t+1}
\geq\Psi'\!\left(\mathfrak{s}_t \right)(\mathfrak{s}_t  - \mathfrak{s}_{t+1}).
\]
Combining the latter with \eqref{eq:noqap-corrected-descent} and \eqref{eq:noqap-key-estimate} yields
\begin{equation}\label{eq:001_dif}
\begin{aligned}
&\frac{\|x^{t+1}-x^t\|^2}{\alpha_{t+1}}
+\alpha_{t+1}\dist^2\bigl(0;\partial\varphi(x^{t+1})\bigr)\\
&\quad\leq\frac{\widehat\Delta_t-\widehat\Delta_{t+1}}{\mathfrak r}
\left[\dist\bigl(0;\partial\varphi(x^t)\bigr)
+\frac1{\psi'\!\left(\sum_{j=t}^{\infty}(\alpha_j^3+\alpha_j\varepsilon_j^2)
\right)}\right].
\end{aligned}
\end{equation}
The construction of $\Gamma_t$, the inequality $(a+b)^2\leq2(a^2+b^2)$,
and \eqref{eq:noqap-step-ratio} ensure that
\begin{equation}\label{eq:002_dif}
    \begin{aligned}
\Gamma_{t+1}^2
&\le 2\alpha_{t+1}
\left(\frac{\|x^{t+1}-x^t\|^2}{\alpha_{t+1}}
+\alpha_{t+1}\dist^2\bigl(0;\partial\varphi(x^{t+1})\bigr)\right)\\
&\le\frac{2\rho}{\mathfrak r}\left(\Gamma_t+\frac{\alpha_t}{\psi'\!\left(\sum_{j=t}^{\infty}\bigl(\alpha_j^3+\alpha_j\varepsilon_j^2\bigr)\right)}\right)
(\widehat\Delta_t-\widehat\Delta_{t+1}).
\end{aligned}
\end{equation}

Consequently, we get the estimates
\begin{align*}
\Gamma_{t+1}
&\le
\left(\Gamma_t+\frac{\alpha_t}{\psi'\!\left(\sum_{j=t}^{\infty}\bigl(\alpha_j^3+\alpha_j\varepsilon_j^2\bigr)\right)}\right)^{1/2}
\left(\frac{2\rho}{\mathfrak r}\left(\widehat\Delta_t-\widehat\Delta_{t+1} \right)\right)^{1/2} \\
&\le\frac{1}{2} \left(\Gamma_t+\frac{\alpha_t}{\psi'\!\left(\sum_{j=t}^{\infty}\bigl(\alpha_j^3+\alpha_j\varepsilon_j^2\bigr)\right)}\right) +  \frac{\rho}{\mathfrak r} \left(\widehat\Delta_t-\widehat\Delta_{t+1}\right).
\end{align*}
Rearranging terms and multiplying by $2$ verifies \eqref{eq:noqap-Gamma-recursion}.\vspace*{0.03in}

{\textbf{Claim~3:}} {\em Let $\mathfrak{s}_t\in]0,\eta[$ and $x^t\in U$ for every $t=t_0,\ldots,T$. Then}
\begin{align}
\sum_{t=t_0}^{T}\Gamma_{t+1}
& \leq\Gamma_{t_0}+\frac{2\rho}{\mathfrak r}\widehat\Delta_{t_0}
+\sum_{t=t_0}^{T}\frac{\alpha_t}{\psi'\!\left(\sum_{j=t}^{\infty}\bigl(\alpha_j^3+\alpha_j\varepsilon_j^2\bigr)\right)}.
\label{eq:noqap-finite-sum}
\end{align}
Indeed, using \eqref{eq:noqap-Gamma-recursion} with subsequent summing and rearranging gives us
\begin{align*}
\sum_{t=k}^{T}\Gamma_{t+1}
&\le\Gamma_k-\Gamma_{T+1}
+\frac{2\rho}{\mathfrak r}
(\widehat\Delta_k-\widehat\Delta_{T+1})
+\sum_{t=t_0}^{T}\frac{\alpha_t}{\psi'\!\left(\sum_{j=t}^{\infty}\bigl(\alpha_j^3+\alpha_j\varepsilon_j^2\bigr)\right)}\nonumber\\
&\leq\Gamma_k+\frac{2\rho}{\mathfrak r}\widehat\Delta_k
+\sum_{t=t_0}^{T}\frac{\alpha_t}{\psi'\!\left(\sum_{j=t}^{\infty}\bigl(\alpha_j^3+\alpha_j\varepsilon_j^2\bigr)\right)},
\end{align*}
which readily verifies the claim.\vspace*{0.03in}

{\textbf{Claim~4:}} {\em The conclusions of the theorem hold}. To complete the proof of the theorem, let $r>0$ be such  that $\mathbb{B}_{r}[\bar x] \subseteq U$. Select $t_0 \in\N$ satisfying
\begin{equation}\label{eq:t0}
\| x^{t_0} - \bar x\| + \Gamma_{t_0}+\frac{2\rho}{\mathfrak r}\widehat\Delta_{t_0}
+\sum_{t=t_0}^{T}\frac{\alpha_t}{\psi'\!\left(\sum_{j=t}^{\infty}\bigl(\alpha_j^3+\alpha_j\varepsilon_j^2\bigr)\right)}< r,
\end{equation}
which exists by \eqref{desing_condi}, \eqref{eq:noqap-small-Gamma}, and the fact that $\bar x$ is an accumulation point of $\{x^t\}$. Without loss of generality, we can assume that $\mathfrak{s}_t \in ]0,\eta[$ for all $t\geq t_0$. In particular, \eqref{eq:t0} implies that $x^{t_0} \in U$. Proceeding by induction, suppose that $x_{t} \in U$ for all $t=t_0, \ldots, T$. It follows from Claim~3 tell us that \eqref{eq:noqap-finite-sum} holds. This tells us that
\begin{align*}
\| x^{T+1} - \bar x\| &\leq \| x^{t_0} - \bar {x}\| + \sum_{t=t_0}^{T} \Gamma_{t+1} \\& \le\| x^{t_0} - \bar x\| + \Gamma_{t_0}+\frac{2\rho}{\mathfrak r}\widehat\Delta_{t_0}+\sum_{t=t_0}^{T}\frac{\alpha_t}{\psi'\!\left(\sum_{j=t}^{\infty}\bigl(\alpha_j^3+\alpha_j\varepsilon_j^2\bigr)\right)} <r.
\end{align*}
In particular, the latter implies that $x^{T+1} \in U$, and thus induction shows that $x^{t} \in U$ for all $t\geq t_0$. Therefore, we get that \eqref{eq:noqap-finite-sum} holds for arbitrary $T\geq t_0$, and consequently $\sum_{t\in \N} \Gamma_t<\infty$. This verifies \eqref{eq:noqap-conclusions} and thus completes the proof.
\end{proof}

\begin{remark}\label{rem:value-attainment} Observe that assumption \eqref{desing_condi} for sufficiently large $p_0$ is needed  to guarantee that the term $\sum_{j=p_0}^\infty(\alpha_j^3+\alpha_j\varepsilon_j^2)$ belongs to the domain of $\psi$. Furthermore, it follows from the proof of the theorem that the continuity of $h$ relative to its domain can be replaced by assuming the existence of
a subsequence $x^{t_j}\to\bar x$ that satisfies the condition
$\varphi(x^{t_j})\to\varphi(\bar x)$ as $j\to\infty$.
\end{remark}

\begin{remark}\label{rem:smooth-no-step-ratio}
If $h\equiv0$, then the adjacent stepsize condition
\eqref{eq:noqap-step-ratio} can be omitted from
Theorem~\ref{Theo_main_convergencePLK} with the fulfillment of all its conclusions. Indeed, the proximal step reduces in this case to the equalities
\[
x^t-x^{t-1}=-\alpha_t\mathtt N
\left(\nabla f(x^{t-1})+e_t\right),
\qquad
e_t:=\frac{1}{\mathtt N}\sum_{i=1}^{\mathtt N}g_i^t
-\nabla f(x^{t-1}).
\]
Combining this with the smoothness of $f$ and estimate \eqref{eq:direct-epoch-discrepancy} gives us
\begin{equation}\label{eq:smooth-corrected-descent-remark}
\begin{aligned}
f(x^t)\le f(x^{t-1})
&-\widetilde{\mathfrak r}\left(
\frac{\|x^t-x^{t-1}\|^2}{\alpha_t}+\alpha_t\|\nabla f(x^{t-1})\|^2\right) +\widetilde n_0
\left(\alpha_t^3+\alpha_t\varepsilon_t^2\right)
\end{aligned}
\end{equation}
for all large $t$ with some constants $\widetilde{\mathfrak r}>0$ and $\widetilde n_0>0$.
The important point here is that the gradient appearing in this inequality is $\nabla f(x^{t-1})$ while not $\nabla f(x^t)$, with the same stepsize $\alpha_t$ that controls the increment from $x^{t-1}$ to $x^t$; cf.\ \eqref{eq:noqap-corrected-descent}).   Therefore, mimicking the above proof, we can proceed in this case without imposing \eqref{eq:noqap-step-ratio}.
\end{remark}

The next theorem provides {\em convergence rates} for Algorithm~\ref{alg:1} under the {\em exponent} version of PLK condition. Observe the elementary power inequality used in the proof.
\begin{equation}\label{eq_0001}
(a+b)^\theta\leq a^\theta+b^\theta,\quad(a,b\ge 0,\quad 0\le\theta<1.
\end{equation}

\begin{theorem}\label{main_rates} Let $\{x^t\}$ be a sequence of iterates generated by Algorithm~{\rm\ref{alg:1}} with
\[
\alpha_t=\frac{a}{(t+b)^\gamma},\qquad
\varepsilon_t=\mathcal O\bigl((t+b)^{-\gamma}\bigr),
\qquad a>0,\quad b\geq 0,\quad\frac12<\gamma\le 1.
\]
In addition to Assumption~{\rm\ref{Assumpt01}}, suppose that $\bar x$ is an accumulation point satisfying the $($exponent$)$ PLK condition with $\psi(s)=Ms^{1-\theta}$, where $M>0$ and $\theta\in]0,1[$. Then we have $x^t\to\bar x$, $\varphi(x^t)\to\varphi(\bar x)$ as $t\to\infty$ with the following rate estimates:
\begin{equation}\label{ratevalues}
\left|\varphi(x^t)-\varphi(\bar x)\right|
=\begin{cases}
\mathcal O\bigl(t^{-(3\gamma-1)}\bigr),
&\frac12<\gamma<1,\quad
0<\theta\leq\dfrac{\gamma}{3\gamma-1},\\
\mathcal O\bigl(t^{-(1-\gamma)/(2\theta-1)}\bigr),
&\frac12<\gamma<1,\quad
\dfrac{\gamma}{3\gamma-1}<\theta<1,\\
\mathcal O(t^{-2}),
&\gamma=1,\quad0<\theta<\frac12,\\
\mathcal O(t^{-2}&\gamma=1,\quad\theta=\frac12,\quad \frac{a \mathtt{N} }{M^2}>1,\\
\mathcal O^-\bigl(t^{- \frac{2a \mathtt{N} }{M^2}} \bigr),
&\gamma=1,\quad\theta=\frac12,\quad0<  \frac{a \mathtt{N} }{M^2}\leq1,\\
\mathcal O\bigl([\ln(t+1)]^{-1/(2\theta-1)}\bigr),
&\gamma=1,\quad\frac12<\theta<1,
\end{cases}
\end{equation}
\begin{equation}\label{ratespoints}
\|x^t-\bar x\|=\begin{cases}
\mathcal O\bigl(t^{-(2\gamma-1)}\bigr),
&\frac12<\gamma<1,\quad
0<\theta\leq\dfrac{\gamma}{3\gamma-1},\\
\mathcal O\bigl(t^{-(1-\gamma)(1-\theta)/(2\theta-1)}\bigr),
&\frac12<\gamma<1,\quad
\dfrac{\gamma}{3\gamma-1}<\theta<1,\\
\mathcal O(t^{-1}),
&\gamma=1,\quad0<\theta<\frac12,\\
\mathcal O(t^{-1}),
&\gamma=1,\quad\theta=\frac12,\quad\frac{a \mathtt{N} }{M^2}>1,\\
\mathcal O^{-}\bigl(t^{-\frac{a \mathtt{N} }{M^2}}\bigr),
&\gamma=1,\quad\theta=\frac12,\quad0< \frac{a \mathtt{N} }{M^2}\le 1,\\
\mathcal O\bigl([\ln(t+1)]^{-(1-\theta)/(2\theta-1)}\bigr),
&\gamma=1,\quad\frac12<\theta<1.
\end{cases}
\end{equation}
\end{theorem}
\begin{proof}
First we verify the convergence $x^t\to\bar x$ and $\varphi(x^t)\to\varphi(\bar x)$ as $t\to\infty$. In our setting, the conditions in \eqref{eq01:prop:stat2} and \eqref{eq:noqap-step-ratio} are satisfied. Choose $\vartheta\in[\theta,1[$ so close to 1 that $\gamma+\vartheta(3\gamma-1)>1$; this is possible due to $\gamma >1/2$. On a sufficiently small positive value interval, we have that the function
\begin{equation}\label{def:hatpsi}
\psi_\vartheta(r)=\widehat M r^{1-\vartheta}\quad\mbox{ with }\;
\widehat M=\frac{M(1-\theta)}{1-\vartheta}
\end{equation}
is another desingularizer, and it is easy to see that $\varphi$ satisfies the basic PLK condition from Definition~\ref{def_PKL} with the desingularizing function $\psi_\vartheta$. Moreover,  Lemma~\ref{Lemma_3.7rate} gives us \eqref{desing_condi} for $\psi_\vartheta$, and hence
Theorem~\ref{Theo_main_convergencePLK} ensures that $x^t\to\bar x$ and $\varphi(x^t) \to \varphi(\bar x)$.

Fix further any $\delta\in]0,1/2[$ and set
\begin{equation}\label{kappa_Theo}
s:=3\gamma-1,\quad\kappa:=\frac12,
\qquad
\lambda_\delta:=\frac{4(\kappa-\delta)a\mathtt N}{M^2},
\end{equation}

which allows us to deduce from		Lemma~\ref{Lemma_3.7rate} that $\mathtt u_t=\mathcal O((t+b)^{-s})$.

By \eqref{eq:noqap-corrected-descent}, the nonnegative quantity $\varphi(x^t)-\varphi(\bar x)+\mathtt u_t$ is eventually nonincreasing and converges to zero. If it vanishes at some sufficiently large index, then the corrected descent inequality yields $x^t=\bar x$ from that index onward, and all the  claimed rates are readily satisfied. We henceforth assume that
\[
\varphi(x^t)-\varphi(\bar x)+\mathtt u_t>0
\]
for all large $t\in\N$. We split the rest of the proof into the following four claims.\vspace*{0.03in}

{\textbf{Claim~1:}} {\em Given any $\delta \in ]0,1/2[$ and $\vartheta\in[\theta,1[$, there exists $t_0 \in \N$ such that
\begin{equation}\label{ineq_claim01:rates}
\begin{aligned}
&\frac{(\kappa-\delta)\mathtt N\alpha_t}{2M^2(1-\theta)^2}
\bigl(\varphi(x^t)-\varphi(\bar x)+\mathtt u_t\bigr)^{2\vartheta}\\
&\quad\leq
\bigl(\varphi(x^{t-1})-\varphi(\bar x)+\mathtt u_{t-1}\bigr)
-\bigl(\varphi(x^t)-\varphi(\bar x)+\mathtt u_t\bigr)\\
&\qquad\quad+
\frac{(\kappa-\delta)\mathtt N\alpha_t}{M^2(1-\theta)^2}
\mathtt u_t^{2\vartheta}\quad\mbox{ for all }\quad t\geq t_0.
\end{aligned}
\end{equation}
If $\vartheta=1/2$, then the sharper estimate
\begin{equation}\label{eq:half-exponent-recursion}
\begin{aligned}
&\left(1+\frac{4(\kappa-\delta)\mathtt N\alpha_t}{M^2}\right)
\bigl(\varphi(x^t)-\varphi(\bar x)+\mathtt u_t\bigr)\\
&\qquad\leq\varphi(x^{t-1})-\varphi(\bar x)+\mathtt u_{t-1}
+\frac{4(\kappa-\delta)\mathtt N\alpha_t}{M^2}\mathtt u_t
\end{aligned}
\end{equation}
also holds for all sufficiently large $t$}.

To verify the claim, fix $\vartheta\in[\theta,1[$ and consider $\psi_\vartheta$ from \eqref{def:hatpsi}. The iterates eventually belong to the PLK neighborhood, and their objective values converge to $\varphi(\bar x)$. If $\varphi(x^t)>\varphi(\bar x)$, then \eqref{eq_0001} and the PLK condition yields
\begin{align*}
\bigl(\varphi(x^t)-\varphi(\bar x)+\mathtt u_t\bigr)^\vartheta & \le
\bigl(\varphi(x^t)-\varphi(\bar x)\bigr)^\vartheta
+\mathtt u_t^\vartheta  \leq \widehat{M}(1-\vartheta)
\dist\bigl(0;\partial\varphi(x^t)\bigr)+\mathtt u_t^\vartheta.
\end{align*}
On the other hand, for $\varphi(x^t)\leq\varphi(\bar x)$ we get $0<\varphi(x^t)-\varphi(\bar x)+\mathtt u_t\leq\mathtt u_t$, which verifies the same bound. Squaring and using $(u+v)^2\le2(u^2+v^2)$ give us
\begin{align*}
&\bigl(\varphi(x^t)-\varphi(\bar x)+\mathtt u_t\bigr)^{2\vartheta} \le
2 \widehat{M}^2(1-\vartheta)^2\dist^2\bigl(0;\partial\varphi(x^t)\bigr)
+2\mathtt u_t^{2\vartheta}.
\end{align*}
It follows from Theorem~\ref{prop:01} and the construction of $\mathtt u_t$ that
\[
\begin{aligned}
(\kappa-\delta)\mathtt N\alpha_t
\dist^2\bigl(0;\partial\varphi(x^t)\bigr)
&\le\varphi(x^{t-1})-\varphi(x^t)
+\mathtt u_{t-1}-\mathtt u_t.
\end{aligned}
\]
The last two estimates verifies \eqref{ineq_claim01:rates}. Since $\alpha_t=a(t+b)^{-\gamma}$ and $\mathtt u_t=\mathcal O((t+b)^{-s})$, this is a recursion covered by Lemma~\ref{lem:rate_sequence} with $q=2\vartheta$ and $\mu=s$ whenever $2\vartheta\ge 1$.

Consider now the case where $\theta=1/2$. If $\varphi(x^t)-\varphi(\bar x)>0$, then the PLK condition tells us that
$\varphi(x^t)-\varphi(\bar x)\le\frac{M^2}{4}\dist^2\bigl(0;\partial\varphi(x^t)\bigr)$,
which yields
\[
\varphi(x^t)-\varphi(\bar x)+\mathtt u_t
\le\frac{\widehat{M}^2}{4}\dist^2\bigl(0;\partial\varphi(x^t)\bigr) +  \mathtt{u}_t.
\]
If $ \varphi(x^t)-\varphi(\bar x)\le 0$, we simply have
$$
0 < \varphi(x^t)-\varphi(\bar x) +\mathtt{u}_t \leq \mathtt{u}_t \leq \frac{\widehat{M}^2}{4}\dist^2\bigl(0;\partial\varphi(x^t)\bigr) +  \mathtt{u}_t.
$$
Thus in both cases for $\theta=1/2$, the estimate in \eqref{eq:half-exponent-recursion} is verified, which therefore completes the proof in Claim~1.\vspace*{0.03in}

{\textbf{Claim~2:}} {\em For every $\vartheta\in[\theta,1[$ satisfying $\gamma+\vartheta s>1$, there exist $C>0$ and $t_0 \in \N$  ensuring the estimate
\begin{equation}\label{eq:dist_bound_claim4}
\begin{aligned}
\|x^t-\bar x\|
&\le C\sqrt{\alpha_t
\bigl(\varphi(x^{t-1})-\varphi(\bar x)+\mathtt u_{t-1}\bigr)}\\
&\quad+C\bigl(\varphi(x^t)-\varphi(\bar x)+\mathtt u_t\bigr)^{1-\vartheta}
+C(t+b)^{-(\gamma+\vartheta s-1)}
\end{aligned}
\end{equation}
for all $t\geq t_0$. Moreover, we have}
\begin{equation}\label{eq:absolute-objective-bound}
\bigl|\varphi(x^t)-\varphi(\bar x)\bigr|
\le\varphi(x^t)-\varphi(\bar x)+\mathtt u_t+\mathtt u_t.
\end{equation}

To verify the claim, we follow the proof of Theorem~\ref{Theo_main_convergencePLK}. Applying \eqref{eq:noqap-finite-sum} with the function
$\psi_\vartheta$ used in Claim~1 gives us the inequalities
\begin{align*}
\|x^t-\bar x\|
&\le\sum_{j=t}^\infty\|x^{j+1}-x^j\|\le\Gamma_t
+C_1\bigl(\varphi(x^t)-\varphi(\bar x)+\mathtt u_t\bigr)^{1-\vartheta}\\
&\quad+C_2\sum_{j=t}^\infty\alpha_j
\left[\sum_{i=j}^\infty(\alpha_i^3+\alpha_i\varepsilon_i^2)\right]^\vartheta\\
& \le\Gamma_t
+C \bigl(\varphi(x^t)-\varphi(\bar x)+\mathtt u_t\bigr)^{1-\vartheta} + C(t+b)^{-(\gamma+\vartheta s-1)}
\end{align*}
for all large $t$ with some $C,C_1,C_2 >0$, where $\Gamma_t$ is defined in \eqref{eq:noqap-small-Gamma}, and where for the last inequality we use Lemma~\ref{Lemma_3.7rate}. It follows from \eqref{eq:noqap-corrected-descent} and \eqref{eq:noqap-value-limit} that
\begin{align*}
\Gamma_t^2
&\le 2\alpha_t\left[\frac{\|x^t-x^{t-1}\|^2}{\alpha_t}
+\alpha_t\dist^2\left( 0;\partial\varphi(x^t)\right) \right]\\
&\leq\frac{2\alpha_t}{\mathfrak r}
\left[\varphi(x^{t-1})	+\mathtt u_{t-1}-(\varphi(x^t)
+\mathtt u_t)\right]
\le\frac{2\alpha_t}{\mathfrak r}
\left[(\varphi(x^{t-1})-\varphi(\bar x)+\mathtt u_{t-1}\right].
\end{align*}
Combining these inequalities and enlarging $C>0$ if needed, we get \eqref{eq:dist_bound_claim4}.
Finally, the estimate in \eqref{eq:absolute-objective-bound} follows from the triangle inequality.\vspace*{0.03in}

{\textbf{Claim~3:}} {\em The results in \eqref{ratevalues} and \eqref{ratespoints} fulfill for $1/2<\gamma<1$}.

To justify this, set $\theta_0:=\frac{\gamma}{3\gamma-1}\in]1/2,1[$ and suppose first that $\theta\leq\theta_0$. Choose $\vartheta=\theta_0$ satisfying the conditions
\[
2\vartheta>1,\qquad
\frac{1-\gamma}{2\vartheta-1}=s,\qquad
\gamma+\vartheta s-1=s(1-\vartheta)=2\gamma-1>0.
\]
Applying Lemma~\ref{lem:rate_sequence} with $q=2\vartheta>1$ and $\mu=s$ to the recursion in Claim~1 yields
\[
\varphi(x^t)-\varphi(\bar x)+\mathtt u_t=\mathcal O(t^{-s}).
\]
Since $\mathtt u_t=\mathcal O(t^{-s})$, the first objective rate in \eqref{ratevalues} follows from \eqref{eq:absolute-objective-bound}. Observe that in \eqref{eq:dist_bound_claim4}, the three terms are bounded, respectively, by the constant multiples of
\[
t^{-(\gamma+s)/2},\qquad
t^{-s(1-\vartheta)},\qquad
t^{-(\gamma+\vartheta s-1)}.
\]
Since $(\gamma+s)/2=2\gamma-1/2>2\gamma-1$ while
$s(1-\vartheta)=\gamma+\vartheta s-1=2\gamma-1$,
Claim~2 gives us the rate $\|x^t-\bar x\|=\mathcal O(t^{-(2\gamma-1)})$.

Now suppose that $\theta>\theta_0$ and choose $\vartheta=\theta$. Put $p:=(1-\gamma)/(2\theta-1)$ with $0<p<s$. Applying Lemma~\ref{lem:rate_sequence} with $q=2\theta>1$ and $\mu=s$ to \eqref{ineq_claim01:rates} brings us to
\[
\varphi(x^t)-\varphi(\bar x)+\mathtt u_t=\mathcal O(t^{-p}).
\]
Since $\mathtt u_t=\mathcal O(t^{-s})$ and $p<s$, we deduce from \eqref{eq:absolute-objective-bound} the second objective rate in \eqref{ratevalues}.
For the iterate estimate, observe that
\begin{align*}
\gamma+\theta s-1-p(1-\theta)
&=\frac{2\theta(\theta s-\gamma)}{2\theta-1}>0,\\
\frac{\gamma+p}{2}-p(1-\theta)&=\frac12>0.
\end{align*}
In particular, we get $\gamma+\theta s>1$, and so Claim~2 applies in this setting. The resulting bounds for the first and third terms in
\eqref{eq:dist_bound_claim4} decay faster than
$t^{-p(1-\theta)}$. Therefore,
$\|x^t-\bar x\|=\mathcal O(t^{-p(1-\theta)})$, which verifies Claim~3.\vspace*{0.03in}

{\textbf{Claim~4:}} {\em If $\gamma=1$, the estimates in \eqref{ratevalues} and \eqref{ratespoints} fulfill}.

Indeed, in this case we have $s=2$ and $\mathtt u_t=\mathcal O(t^{-2})$. Suppose first that $\theta=1/2$ and observe from Claim~1 that
\eqref{eq:half-exponent-recursion} is a particular case of Lemma~\ref{lem:rate_sequence} with $q=1$, $\mu=2$, and $\lambda_\delta=4(\kappa-\delta)a\mathtt N/M^2$. Consequently,
\[
\varphi(x^t)-\varphi(\bar x)+\mathtt u_t
=\begin{cases}
\mathcal O(t^{-2}),&\lambda_\delta>2,\\
\mathcal O(t^{-\lambda_\delta}),&0<\lambda_\delta<2.
\end{cases}
\]
If $\frac{a\mathtt{N}}{M^2} >1$, then we can take $\delta\in]0,\frac{1}{2}[$ such that $\lambda_\delta>2$. On the other hand, for
$\lambda_\delta< \frac{2a\mathtt{N}}{M^2}<2$, the number $\lambda_\delta$ could be arbitrarily close to  $\frac{2a\mathtt{N}}{M^2}$. Then \eqref{eq:absolute-objective-bound} gives us the claimed rates
bounds. When $\vartheta=1/2$, Claim~2 yields
\begin{align*}
\|x^t-\bar x\|
&\le C\sqrt{\alpha_t
\bigl(\varphi(x^{t-1})-\varphi(\bar x)+\mathtt u_{t-1}\bigr)} +C\sqrt{\varphi(x^t)-\varphi(\bar x)+\mathtt u_t}+Ct^{-1}.
\end{align*}
This provides the $\mathcal O(t^{-1})$ rate when $\lambda_\delta>2$ and the
$\mathcal O(t^{-\lambda_\delta/2})$ one when $0<\lambda_\delta<2$. The combined expressions for $0<\lambda_\delta\le2$ in \eqref{ratevalues} and \eqref{ratespoints} reflect both bounds.

Next we focus on the case where $0<\theta<1/2$. Choose $M_0$ so small that $4(\kappa-\delta)a\mathtt N/M_0^2>2$. By shrinking $U$ and $\eta>0$ if needed, we see that the function $r\mapsto M_0\sqrt r$ is also a PLK desingularizer for $\varphi$ at $\bar x$.  Hence the above arguments for $\theta=1/2$, being  applied with $M_0$, give us the objective rate $\mathcal O(t^{-2})$ and the iterate rate $\mathcal O(t^{-1})$.

Finally, we examine the case where $1/2<\theta<1$. Take $\vartheta=\theta$ in Claim~1 and apply Lemma~\ref{lem:rate_sequence} with $q=2\theta>1$, $\mu=2$, and $\gamma=1$ to obtain
\[
\varphi(x^t)-\varphi(\bar x)+\mathtt u_t
=\mathcal O\bigl([\ln(t+1)]^{-1/(2\theta-1)}\bigr).
\]
Since $\mathtt u_t=\mathcal O(t^{-2})$, Claim~2 readily gives the last objective rate in \eqref{ratevalues}. For the iterate estimate, use \eqref{eq:dist_bound_claim4} with
$\vartheta=\theta$. The first term is bounded by
\[
t^{-1/2}[\ln(t+1)]^{-1/[2(2\theta-1)]},
\]
and the third term is of order $t^{-2\theta}$.
Both decay faster than the bound
$[\ln(t+1)]^{-(1-\theta)/(2\theta-1)}$ for the second term. Therefore,
\[
\|x^t-\bar x\|=\mathcal O\bigl([\ln(t+1)]^{-(1-\theta)/(2\theta-1)}\bigr),
\]
which verifies Claim~4 and thus completes the proof of the theorem.
\end{proof}

\begin{remark}\label{Remark:weaklyconvex}
When $h$ is $\rho$-weakly convex, Remark~\ref{rem:weakly-convex-descent} provides the improved coefficient $\mathfrak b_t=(1-\delta)\mathtt N\alpha_t$. Following the proof of Theorem~\ref{main_rates} (with $\kappa=1$ in \eqref{kappa_Theo}), we can consider $\lambda_\delta:=\frac{4a(1-\delta)\mathtt N}{M^2}$. This yields the improved sufficient threshold $\frac{a\mathtt N}{M^2}>\frac12$ for the convergence rates, i.e., we get for $\gamma =1$ and $\theta=1/2$ that
\begin{align*}
\left|\varphi(x^t)-\varphi(\bar x)\right|
&=\begin{cases}\mathcal O(t^{-2}),
&\gamma=1,\quad\theta=\frac12,\quad \frac{a \mathtt{N} }{M^2}>\frac{1}{2},\\
\mathcal O\bigl(t^{-\lambda} \bigr),
&\gamma=1,\quad\theta=\frac12,\quad0<\frac{\lambda}{4}<\frac{a \mathtt{N} }{M^2}\leq\frac{1}{2}.
\end{cases} \\
\|x^t-\bar x\|&=\begin{cases} 	\mathcal O(t^{-1}),
&\gamma=1,\quad\theta=\frac12,\quad \frac{a \mathtt{N} }{M^2}>\frac{1}{2},\\
\mathcal O\bigl(t^{-\lambda/2}\bigr),
&\gamma=1,\quad\theta=\frac12,\quad0<\frac{\lambda}{4}<\frac{a \mathtt{N} }{M^2}\le\frac{1}{2}.
\end{cases}
\end{align*}
\end{remark}

Our final result establishes improved convergence rates for the objective values when the accumulation  point is a local minimizer.

\begin{theorem}
\label{thm:improved-local-rates} In addition to the assumptions of Theorem~{\rm\ref{main_rates}} with $0<\theta\leq1/2$, $1/2<\gamma\leq1$, and
$(\theta,\gamma)\neq (1/2,1)$, suppose that
$\bar x$ is a local minimizer of $\varphi$.
Then the objective gap is eventually nonnegative and satisfies the condition
\begin{equation}\label{eq:local-minimum-improved-rates}
\varphi(x^t)-\varphi(\bar x)=\mathcal O(t^{-\gamma/\theta})\;\mbox{ as }\;t\to\infty.
\end{equation}
\end{theorem}
\begin{proof}
Theorem~\ref{main_rates} gives us $x^t\to\bar x$. Since $\bar x$ is a local minimizer, we have
$\varphi(x^t)-\varphi(\bar x)\ge 0$ for all large $t$. If $\varphi(x^t)>\varphi(\bar x)$, the PLK condition yields
\[
\dist^2\big(0;\partial\varphi(x^t)\big)
\ge\frac{[\varphi(x^t)-\varphi(\bar x)]^{2\theta}}
{M^2(1-\theta)^2},
\]
where the inequality holds trivially if $\varphi(x^t) =\varphi(\bar x)$. It follows from Theorem~\ref{prop:01} that there exists $D>0$ such that
\begin{equation}\label{eq:raw-local-rate-recursion}
\begin{aligned}
&\varphi(x^t)-\varphi(\bar x)+
\frac{(\frac{1}{2}-\delta)a\mathtt N}
{M^2(1-\theta)^2(t+b)^\gamma}
[\varphi(x^t)-\varphi(\bar x)]^{2\theta}\\
&\le
\varphi(x^{t-1})-\varphi(\bar x)
+\frac{D}{(t+b)^{3\gamma}}
\end{aligned}
\end{equation}
for large $t\in\N$. Since $q=2\theta\in]0,1]$ and $(q,\gamma)\neq(1,1)$,
Lemma~\ref{lem:local-rate-scalar} applies with
$\lambda =(\frac{1}{2}-\delta)a\mathtt N/[M^2(1-\theta)^2]$. This verifies \eqref{eq:local-minimum-improved-rates} and thus completes the proof.
\end{proof}

\begin{remark} The obtained
Theorem~\ref{main_rates} reveals a {\em trade-off} between the decay of the stepsizes and the PLK exponent $\theta$. Table~\ref{tab:optimal-rates} summarizes the choices of $\gamma$ that yield the best polynomial convergence rates
for $\theta\neq 1/2$. The critical case where $\theta=1/2$ with $\gamma=1$ deserves
particular attention. In this setting, the convergence
rates depend on the multiplicative constant through $\frac{4 \kappa a\mathtt N}{M^2}$. Observe also that
the sharper descent estimate \eqref{eq:sharp-two-term}
allows us to use $\kappa=1/2$ for  general regularizers $h$ and $\kappa=1$ for weakly convex ones. Since $\delta>0$ may be arbitrarily small, this tells us that a sufficient threshold for $\lambda_\delta>2$
becomes $a\mathtt N>M^2$ in the general case and
$a\mathtt N>M^2/2$ in the weakly convex case. Neither threshold is claimed to be a necessary one.

Note finally that when the PLK constant $M$ is unknown, the parameter $a$ in the stepsize may not be sufficiently large to guarantee the fastest convergence rate
\end{remark}\vspace*{-0.1in}

\begin{table}[htbp]
\centering
\caption{Optimal polynomial stepsize exponents and convergence rates under
the PLK condition, with and without local minimality of $\bar x$.}\label{tab:optimal-rates}
\renewcommand{\arraystretch}{1.7}
\setlength{\tabcolsep}{3pt}
\footnotesize
\begin{tabular}{@{}ccccc@{}}
\toprule
\textbf{PLK exponent} & \textbf{Optimal stepsize}
& \multicolumn{2}{c}{$|\varphi(x^t)-\varphi(\bar x)|$}
& \textbf{Iterate rate} \\
\cmidrule(lr){3-4}
$\theta$ & $\gamma_{\mathrm{opt}}$
& \textbf{General PLK} & \textbf{Local minimizer}
& $\|x^t-\bar x\|$ \\
\midrule
$\displaystyle 0<\theta<\frac12$
& $1$
& $\mathcal O(t^{-2})$
& $\displaystyle\mathcal O\bigl(t^{-1/\theta}\bigr)$
& $\mathcal O(t^{-1})$ \\[4pt]
$\displaystyle\frac12<\theta<1$
& $\displaystyle\frac{\theta}{3\theta-1}$
& $\displaystyle\mathcal O\bigl(t^{-1/(3\theta-1)}\bigr)$
& $\displaystyle\mathcal O\bigl(t^{-1/(3\theta-1)}\bigr)$
& $\displaystyle\mathcal O\bigl(t^{-(1-\theta)/(3\theta-1)}\bigr)$ \\
\bottomrule
\end{tabular}
\par\smallskip
\begin{minipage}{0.98\linewidth}
\footnotesize
\emph{Note.} The local minimizer objective rate for
$0<\theta<1/2$ follows from
Theorem~\ref{thm:improved-local-rates}.
All other rates follow from Theorem~\ref{main_rates};
the iterate-rate column applies to both cases.
At $\theta=1/2$, the rates $\mathcal O(t^{-2})$
and $\mathcal O(t^{-1})$ hold for $\gamma=1$
when the corresponding harmonic coefficient satisfies $\lambda_\delta>2$.
\end{minipage}
\end{table}\vspace*{-0.1in}

\section{Conclusions and Future Research}\label{sec:7}
This paper analyzes a novel proximal inexact random reshuffling method for the composite problem~\eqref{Problem:P} in which a finite sum of
components with Lipschitz continuous gradients is combined with an l.s.c., prox-bounded, and possibly nonconvex regularizer. Under Assumption~\ref{Assumpt01}, we established subsequential stationarity, convergence of the whole sequence under the basic PLK condition, and convergence rates under the exponent one. All statements are pathwise in the sense that they hold for every sequence of permutations and every sequence of inexact gradients within the prescribed tolerances, with no assumption on the distribution of the
reshuffling.

The results rest on two estimates. The within-epoch discrepancy bound controls the epoch error and the gradient tolerance together, and the optimality condition for the proximal step relates the stationarity residual to the increment between consecutive epochs. Being combined, they give us an approximate descent inequality with explicit coefficients whose multiplicative and additive perturbations are of order $\alpha_t^3$ and
$\alpha_t^3+\alpha_t\varepsilon_t^2$. When these perturbations are summable, the objective values converge, although not necessarily monotonically, the weighted increments and stationarity residuals form
convergent series, and every accumulation point is stationary (Theorem~\ref{thm:cluster-stationarity}); none of this requires either the PLK condition, or bounded iterates. Under the PLK condition at a single
accumulation point, together with the adjacent stepsize ratio condition and the weighted error tail condition~\eqref{desing_condi}, the sequence
has finite length and converges to that point
(Theorem~\ref{Theo_main_convergencePLK}). The rescaling Lemma~\ref{lem:noqap-rescaling} removes the quasi-additivity assumption on the desingularizing function $\psi$ without replacing it by another
condition on $\psi$, and the verification, carried out on the corrected quantity $\varphi(x^t)+\mathtt u_t$, does not require the objective values to approach their limit from above. No a priori bound on the trajectory is imposed. Level-boundedness is one sufficient way to guarantee that an accumulation point exists, the continuity assumption
on the regularizer $h$ relative to its domain can be weakened as described in Remark~\ref{rem:value-attainment}, and the ratio condition is unnecessary when $h\equiv0$.

For polynomial stepsizes with exponent $\gamma\in]\frac12,1]$,
tolerances $\varepsilon_t=\mathcal O(\alpha_t)$, and a desingularizing
function $\psi(s)=Ms^{1-\theta}$, Theorem~\ref{main_rates} provides
rates for both the objective values and the iterates in every regime of $(\gamma,\theta)$ summarized in Table~\ref{tab:optimal-rates}. The harmonic stepsize with $\theta=\frac12$ is the critical case, where the rates depend on $a\mathtt N/M^2$; the sharpened descent estimate yields the sufficient thresholds $a\mathtt N>M^2$ in general and $a\mathtt N>M^2/2$ for weakly convex regularizers, neither of which is
claimed to be necessary. When the limit is a local minimizer, a comparison based on the raw objective gap, which no longer involves the
perturbation tail $\mathtt u_t$, gives us the faster objective-value rate of Theorem~\ref{thm:improved-local-rates} while the iterate rates remain
unchanged.

The scope of these conclusions is delimited by the pathwise nature of
the analysis: since no independence or law of the permutations is
exploited, the estimates do not quantify any statistical improvement
caused by randomization itself. The thresholds for the fastest harmonic rates involve the PLK constant $M$; when it is unknown, the stepsize
parameter $a$ may fall below them, and only the slower polynomial rates of Theorem~\ref{main_rates} are guaranteed. Quantitative estimates and
complexity bounds that account for the cost of computing inexact component directions as well as a numerical evaluation of the trade-offs identified here, are natural directions for future work.

\section*{Declarations}

\paragraph{Funding} Research of the first author was partly supported by the US National Science Foundation under grant DMS-2204519 and by the Australian Research Council under Discovery Projects DP-250101112. Research of the second author was partly supported by ANID (Chile) through Fondecyt Regular grants No.~1240120, No.~1261728 and No.~1240335, CMM BASAL funding for the Center of Excellence FB210005, and the project ECOS230027.

\paragraph{Conflict of interest} The authors declare no competing interests.

\appendix
\section{Proof of Lemma~\ref{lem:rate_sequence}}\label{AppendixA}

\begin{proof}
Without loss of generality, we can assume that the recursion inequality \eqref{eq:implicit-scalar-recursion} holds for all $t\in \N$. Considering first the case where $0<\gamma<1$. We claim that there exists $C>0$ and $t_0 \in \N$ such that  $s_{t_0}\leq B_{t_0}$ and
\begin{equation}\label{eq:comparison-supersolution}
B_t+\frac{\lambda}{(t+b)^\gamma}B_t^q
\ge
B_{t-1}+\frac{D}{(t+b)^{\gamma+q\mu}}\;\text{ for all }\;t>t_0,
\end{equation}
where $B_t=C(t+b)^{-p}$. Indeed, let $p: = \min\left\{\mu,\frac{1-\gamma}{q-1}\right\}$ (with the convention $\frac{1-\gamma}{0}=\infty$). By the mean value theorem, there exist $c_p>0$ and
$t_0\in\mathbb N$ such that
\begin{equation}\label{eq:comparison-difference}
(t-1+b)^{-p}  -(t+b)^{-p}
\le\frac{c_p}{(t+b)^{p+1}}
\qquad\text{for all }\;t\ge t_0+1.
\end{equation}
If $q=1$, we get $p=\mu$, and $\gamma<1$ allows us to
increase $t_0$ so that
\begin{equation}\label{eq:lamgamma}
c_\mu(t+b)^{\gamma-1}
\le\frac{\lambda}{2}
\qquad\text{for all }\;t\ge t_0+1.
\end{equation}
Choosing $C :=
\max\left\{\frac{2D}{\lambda},
s_{t_0}(t_0+b)^\mu\right\}$ gives us $s_{t_0}\le B_{t_0}$, and then it follows from \eqref{eq:comparison-difference} and \eqref{eq:lamgamma} that
\[
B_{t-1}-B_t\le\frac{\lambda C}{2(t+b)^{\gamma+\mu}}, \text{ for every }\;t\ge t_0+1.
\]
The above inequality implies that for every
$t\ge t_0+1$, we have the relationships
\[
\begin{aligned}
B_{t-1}-B_t
+\frac{D}{(t+b)^{\gamma+\mu}}\le
\frac{\lambda C/2+D}{(t+b)^{\gamma+\mu}}\le
\frac{\lambda C}{(t+b)^{\gamma+\mu}}=
\frac{\lambda}{(t+b)^\gamma}B_t,
\end{aligned}
\]
which therefore verify \eqref{eq:comparison-supersolution} for $q=1$.

Consider now the case where $q>1$. By the construction of $p$, we get $p\leq\mu$ and
$\gamma+qp\le p+1$. This allows us to choose $C$ so
large that
\[
\lambda C^q\geq c_pC+D,\qquad
C\ge s_{t_0}(t_0+b)^p,
\]
which yield $s_{t_0}\le B_{t_0}$. For $t\geq t_0 +1 $, it follows from \eqref{eq:comparison-difference} and the choice of $C$ that
\[
\begin{aligned}
&B_{t-1}-B_t+\frac{D}{(t+b)^{\gamma+q\mu}} \le
\frac{c_pC}{(t+b)^{p+1}}
+\frac{D}{(t+b)^{\gamma+q\mu}} \le
\frac{c_pC+D}{(t+b)^{\gamma+qp}} \le\frac{\lambda}{(t+b)^\gamma}B_t^q.
\end{aligned}
\]

Next we verify by induction that $s_t\le B_t$ for all $t\geq t_0$. When $t=t_0$, it holds by construction. Suppose that $s_{t-1}\leq B_{t-1}$ for some $t\ge t_0+1$. Using \eqref{eq:implicit-scalar-recursion}
and \eqref{eq:comparison-supersolution}, we obtain the inequalities
\[
\begin{aligned}
s_t+\frac{\lambda}{(t+b)^\gamma}s_t^q
\le
s_{t-1}+\frac{D}{(t+b)^{\gamma+q\mu}}\leq
B_{t-1}+\frac{D}{(t+b)^{\gamma+q\mu}}\leq
B_t+\frac{\lambda}{(t+b)^\gamma}B_t^q.
\end{aligned}
\]
Since the mapping $u\mapsto u+\frac{\lambda}{(t+b)^\gamma}u^q$ is strictly increasing on $[0,\infty[$, this ensures that $s_t\le B_t$,
which completes the proof for $\gamma\in]0,1[$.

Finally, we focus on the case where $\gamma=q=1$. Define $P_t:=
\prod_{j= 1}^{t}
\left(1+\frac{\lambda}{j+b}\right) $ for $t\geq 1$.
Multiplying \eqref{eq:implicit-scalar-recursion}
by $P_{t-1}$ yields
\[
P_ts_t\le P_{t-1}s_{t-1}+\frac{DP_{t-1}}{(t+b)^{1+\mu}}.
\]
Using recursively the above inequality, we obtain
\[
P_ts_t\le P_1s_{1}+D\sum_{j=2}^{t}\frac{P_{j-1}}{(j+b)^{1+\mu}}.
\]
Since $c_1t^{\lambda} \le  P_t  \le c_2 (t+1+b)^{\lambda}$ for some constant $c_1,c_2 >0$, this implies that
\[
c_1t^{\lambda} s_t  \le P_1s_{1} + c_0 D \sum_{j=2}^t (j+b)^{\lambda - \mu -1} = \begin{cases}
\mathcal O(1),&\lambda<\mu,\\
\mathcal O(t^{\lambda-\mu}),&\lambda>\mu,
\end{cases}
\]
and therefore brings us to the estimate
\[
s_t=
\begin{cases}
\mathcal O(t^{-\mu}),&\lambda>\mu,\\
\mathcal O(t^{-\lambda}),&0<\lambda<\mu.
\end{cases}
\]
The remaining case where $\gamma=1$ and $q>1$ is proved similarly to the previous one by considering $B_t=C[\ln(t+b+1)]^{-p}$ with $p=({q-1})^{-1}$.
\end{proof}

\section{Proof of Lemma~\ref{lem:local-rate-scalar}}\label{AppendixB}
\begin{proof}
If $q=1$, the assumptions require that $\gamma<1$.
Apply Lemma~\ref{lem:rate_sequence} with $\mu=2\gamma$, and so $\gamma+q\mu=3\gamma$.
This gives us $s_t=\mathcal O(t^{-2\gamma})$ as claimed. It remains to consider the case where $0<q<1$. Since $3\gamma>1$, the perturbation series is summable, and Lemma~\ref{Prop:Robbins–Siegmund} applied to \eqref{eq:local-scalar-recursion} ensures that that $s_t$ converges and $\sum \lambda s^q_t (t+b)^{-\gamma}<\infty$. By $\gamma\in]1/2,1]$, we have $\sum (t+b)^{-\gamma}=\infty$, which forces $s_t\to 0$ as $t\to\infty$.

Set now $p:=2\gamma/q$. The verification in this case is similar to the proof of Lemma~\ref{lem:rate_sequence}: it consists in considering $B_t=C(t+b)^{-p}$ and adjusting $C>0$ such that there exists $t_0\in\N$ satisfying $s_{t_0}\le B_{t_0}$ and
\begin{equation}\label{eq:lem:local-rate-scalar}
B_{t-1}+\frac{D}{(t+b)^{3\gamma}}
\le B_t+\frac{\lambda}{(t+b)^\gamma}B_t^q
\quad\text{for all }\;t\ge t_0+1.
\end{equation}
To select $C$, we first note that $p+1-3\gamma=2\gamma\left(q^{-1}-1\right)+1-\gamma>0$ and employ the mean value theorem to find $c_p>0$ such that
\[
({t-1}  + b)^{-p}-({t}  + b)^{-p}\le c_p(t+b)^{-p-1}
\]
whenever $t$ is sufficiently large. Further, we fix $c_0>0$ with $\lambda c_0^q\ge 2D$ and taking into account that $0<1-q<1$ by \eqref{eq_0001} and $pq=2\gamma$, arrive at the estimate
\[
(s_{t}(t+b)^p+c_0)^{1-q}(t+b)^{-(p+1-3\gamma)}
\le s_{t}^{1-q}(t+b)^{\gamma-1}+c_0^{1-q}(t+b)^{-(p+1-3\gamma)}
\longrightarrow 0.
\]
This allows us to select $t_0\in\N$ satisfying
\[
c_p(s_{t_0}(t_0+b)^p+c_0)^{1-q}(t_0+b)^{-(p+1-3\gamma)}\leq\frac{\lambda}{2}.
\]
Define $C:=s_{t_0}(t_0+b)^p+c_0$, which yields $s_{t_0}\le B_{t_0}$. For all $t\ge t_0+1$, we get
\[
\begin{aligned}
B_{t-1}-B_t+\frac{D}{(t+b)^{3\gamma}}
&\le c_pC(t+b)^{-p-1}+\frac{D}{(t+b)^{3\gamma}}\\
&\le\frac{\lambda C^q/2+D}{(t+b)^{3\gamma}}
\le\frac{\lambda C^q}{(t+b)^{3\gamma}}
=\frac{\lambda}{(t+b)^\gamma}B_t^q.
\end{aligned}
\]
Therefore, the claimed result now follows by induction, \eqref{eq:lem:local-rate-scalar}, and the monotonicity of the mapping
$u\mapsto u+\lambda(t+b)^{-\gamma}u^q$.
\end{proof}

\bibliographystyle{siam}

\end{document}